\documentclass{amsart}

\usepackage{amssymb}

\title[Applications of modular representation theory]
{Potential applications of modular representation theory to quantum mechanics}
\author{Robert Arnott Wilson}
\address{Queen Mary University of London}
\email{r.a.wilson@qmul.ac.uk}

\date{First draft, 28th May 2021; this version, 26th July 2021.}
\begin{document}
\begin{abstract}
There 
is a unique finite group that lies inside the $2$-dimensional unitary group but not in the special unitary group,
and maps by the symmetric square  
to an irreducible subgroup of the $3$-dimensional real special orthogonal group.
In an earlier paper I showed how 
the representation theory of this group over the real numbers gives rise to much of the 
structure of the standard model of particle physics, but with a number of added twists.
In this theory the group is quantised, but the
representations are not. 
In this paper I consider how a quantisation of the representations might lead to a
more fundamental theory.
\end{abstract}

\maketitle
\tableofcontents
\section{Introduction}
\label{intro}
\subsection{The problem}
\label{problem}
Many, if not all, of the fundamental processes of nature, such as radioactive decay of unstable atomic nuclei, or the transitions of
electrons between ground states and excited states, are discrete. Nevertheless, our best current model of these processes,
the standard model of particle physics \cite{Griffiths}, based on the principles of quantum field theory \cite{Woit}, that generalise the original
theory of quantum mechanics, is an essentially continuous theory. 
This tension between the discrete and continuous bothered Einstein from very early days \cite{Einstein1935}, and still bothers people today,
because it leads to serious problems in interpreting the theory.
In practice, the theory works extremely well, so the default position of modern physicists is to give up worrying about the interpretation,
and simply get on with the job of measuring and calculating. 

But if one wants to attack the deep theoretical problems, like how to
quantise gravity, or how to explain the masses of the elementary particles, or how to explain dark matter and dark energy, then 
it is no longer good enough just to ignore the problem of interpretation. For example, Rovelli interprets the wave-function as
describing relations between particles \cite{RQM}, rather than particles themselves, in order to side-step the measurement problem,
and quantise gravity using loops \cite{Rovelli}.
Other approaches to quantum gravity \cite{verlinde,hossenfelder,mcculloch} are based on
quantisation via information, which similarly requires elementary particles to carry a finite number of bits of information,
rather than a continuous wave-function.

None of these approaches really grapples with the central problem of making the continuous discrete, however. A first step in this
direction might be to replace some of the continuous Lie groups in the standard model by finite groups. In \cite{perspective}
I surveyed a number of finite groups that might be useful in this regard, without coming to any firm conclusions. In \cite{finite}
I investigated one particular finite group in detail, namely the binary tetrahedral group, and found that its real group algebra
contained most of the structures necessary for the standard model of particle physics, 
as well as a possible algebraic basis for a quantum theory of
gravity. 

However, it seemed to me that a more satisfactory theory could be obtained by adjoining an automorphism of order
$2$ to this group, for the purpose of describing the weak interaction. I therefore repeated the exercise \cite{octahedral}
for this extended group. At that stage it became apparent that making the group finite, but keeping the representations
continuous, was only doing half the job. It is therefore my purpose in this paper to make the representations discrete
(or even finite) as well,
by working over the integers, or a finite field (or fields), instead of real or complex numbers.

\subsection{A toy model}
\label{toy}
To illustrate the methods in the simplest possible case, consider the representations of the group $Sym(2)$ of order $2$,
with coefficients in the field $\mathbb F_2$ of order $2$. The elements of the field are $0$ and $1$, subject to all the usual
rules of arithmetic except that $1+1=0$. The group has two elements, $1$ and $d$, say, 
which can be written as
$2\times 2$ matrices
\begin{eqnarray}
1=\begin{pmatrix}1&0\cr0&1\end{pmatrix},&&
d=\begin{pmatrix}0&1\cr1&0\end{pmatrix}.
\end{eqnarray}
These matrices generate a $2$-dimensional algebra, which has a unique $1$-dimensional ideal (that is, $d$-invariant subalgebra), that is 
generated by $1+d$.

Now in \cite{octahedral}, the element $d$ here is used in place of $\gamma_5$ in the standard model, so that $1+d$ is the
appropriate analogue of the projection with $(1-\gamma_5)/2$. But we cannot divide by $2$, since $2=1+1=0$ in $\mathbb F_2$,
so in place of the standard model projection satisfying
\begin{eqnarray}
(1-\gamma_5)^2 &=& 2(1-\gamma_5)
\end{eqnarray}
we have 
\begin{eqnarray}
(1+d)^2&=&0.
\end{eqnarray}
This provides a fundamental mathematical reason why there is no corresponding `right-handed' projection
with $(1+\gamma_5)/2$, and why left-handed spinors and right-handed spinors seem to behave
in completely different ways. The `right-handed' part is glued on top of the `left-handed' part, so cannot exist on its own.

The group algebra over $\mathbb F_2$ is $2$-dimensional, spanned by $1$ and $d$, so has four elements, $0$, $1$, $d$ and $1+d$.
But it is not isomorphic to the sum of two copies of the field, 
that one might (naively) expect by analogy with real or complex
representation theory. 
The `left-handed' part consists of $0$ and $1+d$, and is quite different from the `right-handed'
part in $1$ and $d$. 

One can also transform to a `chiral basis' such that
\begin{eqnarray}
1=\begin{pmatrix}1&0\cr0&1\end{pmatrix},&&
d=\begin{pmatrix}1&0\cr1&1\end{pmatrix},
\end{eqnarray}
acting on row vectors so that $d$ maps $(x,y)$ to $(x+y,y)$.  Then we see more clearly that $d$ acts on the left hand coordinate $x$
but not on the right hand coordinate $y$.

This toy example already shows that modular representation theory has exactly the right properties to model some of the
most puzzling aspects of particle physics, such as the lack of symmetry between left-handed and right-handed spinors.
Of course, we haven't got any spinors in the model yet, so we need to extend to larger groups for complete information.
Let us work our way gradually up the chain of quotient groups through $Sym(3)$ and $Sym(4)$ to the full group $2.Sym(4)$.

\section{A three-point model}
\label{3points}
\subsection{A bigger toy}
\label{bigtoy}
Consider next the group $Sym(3)$ of order $6$. This group is isomorphic to the group $GL(2,2)$ of all invertible $2\times 2$
matrices over the field $\mathbb F_2$ of order $2$. For compatibility of notation with \cite{octahedral} let us take group generators
\begin{eqnarray} 
d=\begin{pmatrix}1&0\cr1&1\end{pmatrix},&&
w=\begin{pmatrix}0&1\cr1&1\end{pmatrix}, 
\end{eqnarray}
so that the other non-identity elements of the group are
\begin{eqnarray} w^2=w^{-1}&=&\begin{pmatrix}1&1\cr1&0\end{pmatrix},\cr
wd=dw^{-1}&=&\begin{pmatrix}1&1\cr0&1\end{pmatrix},\cr
dw=w^{-1}d&=&\begin{pmatrix}0&1\cr1&0\end{pmatrix}. 
\end{eqnarray}

The algebra generated by these matrices is the full $4$-dimensional matrix algebra, and forms a direct summand of the
$6$-dimensional group algebra over $\mathbb F_2$. The other summand is a copy of the $2$-dimensional algebra
described in Section~\ref{toy} above, in which $w$ acts as the identity matrix. These two summands of the group algebra can be
separated by the pair of orthogonal idempotents
\begin{eqnarray}
e:=w+w^2,&&f:=1+e.
\end{eqnarray}
The corresponding irreducible representations are the $2$-dimensional representation we started with, and the trivial 
representation. 

In terms of physical interpretation, we might want the triplet symmetry $w$ to act on the three generations, or the three colours, or perhaps both.
The breaking of the group algebra into $2+4$ could then represent a breaking of generation symmetry into a first generation in the
$2$-dimensional subalgebra, and the other two generations in the matrix algebra. Or it might represent colour confinement, such that the three
colours form a $2$-dimensional representation, and the mixing of colours takes place in the matrix algebra, leaving the other
$2$ dimensions colourless. Or perhaps both, in different contexts. Or perhaps it is just a toy, and has no such interpretations.

\subsection{Representations modulo $3$}
\label{S3mod3}
Now the group $Sym(3)$ has order divisible by $3$, which means that something interesting happens to the representation
theory also in $3$-modular arithmetic. So let us switch attention to the field $\mathbb F_3$ of order $3$, consisting of elements
$0$, $1$ and $-1$ subject to the rule $1+1=-1$.
Here there are again two irreducible representations, both of dimension $1$. One is the trivial representation $1^+$, and the other
is the alternating representation $1^-$, in which $d$ acts as $-1$. 

The group algebra does not split as a direct sum of subalgebras,
but as a representation of $Sym(3)$  it splits as the sum of two indecomposable 
$3$-dimensional representations. One of these may be obtained from the permutation representation of the group, given by matrices
\begin{eqnarray}
w=\begin{pmatrix}0&1&0\cr0&0&1\cr1&0&0\end{pmatrix},&&
d=\begin{pmatrix}1&0&0\cr0&0&1\cr0&1&0\end{pmatrix},
\end{eqnarray}
and the other by changing the sign of $d$. Note that $d$ has determinant $-1$ in the permutation representation,
and determinant $+1$ in the other (so-called monomial) representation.

The abstract structure of this representation becomes clearer if we change to a basis in which $w$ is lower-triangular, and $d$ is diagonal, such as
\begin{eqnarray}
w=\begin{pmatrix}1&0&0\cr1&1&0\cr-1&1&1\end{pmatrix},&&
d=\pm\begin{pmatrix}1&0&0\cr0&-1&0\cr0&0&1\end{pmatrix}.
\end{eqnarray}
It then becomes clear that $w-1$, or equivalently $w-w^2$,
plays a role modulo $3$ that is analogous to the role of $1+d$ modulo $2$, and is therefore likely to
be relevant to the modelling of the strong force. We have
\begin{eqnarray}
(w-w^2)^3&=& 0
\end{eqnarray}
so that the $3$-dimensional algebra has a basis $1, w-w^2, (w-w^2)^2$ that may describe the generation symmetry and/or colour symmetry
in a more revealing manner than the basis $1,w,w^2$. 

In particular, there appears to be a hierarchy of generations, rather than a complete symmetry, and there is an alternating sign as we progress
up or down this hierarchy. Since this sign is determined by $d$, which has something to do with the weak interaction, this suggests a possible role
for a hierarchy of alternating up and down quarks, such as down-up-strange or charm-bottom-top. The former played an important role
in Gell-Mann's original eightfold way \cite{GellMann}, which may thereby acquire a more solid mathematical foundation in the
representation theory I am describing here \cite{Alperin,ABC}.

There is one further mathematical issue that this example throws up, and that is that, unlike the previous examples, the structure of
tensor products of representations is not determined by the character table. That is, 
while it is easy to determine the composition factors of the
representation, 
it is not so easy to determine which order they come in. We have decomposed the group algebra into two so-called
principal (or projective) indecomposable modules (PIMs) of $3$ dimensions each, but there are also two irreducibles $1^\pm$, and two
indecomposables of dimension $2$, obtaining by gluing $1^+$ on top of $1^-$ or vice versa. 

There is just one 
piece of
information 
we are missing, and that is the tensor product of the $2$-dimensional indecomposables.  
The signs are easy to determine, and by direct matrix calculations we find that
\begin{eqnarray}
\begin{matrix} 1^+\cr 1^-\end{matrix} \otimes \begin{matrix}1^+\cr 1^-\end{matrix} &=& 1^- + \begin{matrix}1^+\cr 1^-\cr 1^+\end{matrix}.
\end{eqnarray}
The first component is the anti-symmetric square, and the second is the symmetric square.
Since the latter is a projective module, it gives us no information that isn't already in the character table, so that the tensor
product structure of the representation theory is often expressed with all the projectives omitted. 

I speculate that this particular equation may have something interesting to tell us about what happens when two particles
interact, via a tensor product. However, this is still only a toy model, and we must build the spinors before we can discuss this idea in detail.
Nevertheless, it is the structure of these two PIMs that underlies all the modelling of mass and charge, and the three generations
and three colours, that I propose in this paper. In the next two subsections I summarise in general terms how this might work, while leaving
the detailed discussions of the spinors to Section~\ref{spinors}.

\subsection{Three generations of mass}
\label{3gens}
In the absence of any information about charge, spin, momentum and so on,
the three generations of fermions would seem to be modelled by the permutation representation of $Sym(3)$ on three letters.
As a real representation, this splits as $1^++2^0$, and it is natural to associate the scalar with total energy, and use the representation
$2^0$ as a `mass plane', to distinguish the three generations. Then one can construct a group $SL(2,\mathbb R)$ to act as a gauge
group on the mass plane, as the basis for a theory of the weak interaction.

But this process cannot be carried out in the category of integral representations, because if it could, then it could also be carried out
over every finite field, and we have seen that it cannot.
As a finite representation, over $\mathbb F_3$, the representation is indecomposable, and cannot be split up 
as $1+2$. Its
structure is that of the PIM $P_1^+$, corresponding to the irreducible $1^+$, and with structure $1^+.1^-.1^+$, reading from
bottom to top, say. Corresponding to the `mass plane', we have either the submodule $1^+.1^-$ or the quotient module $1^-.1^+$,
but not both at the same time. 

The mass, as a scalar, corresponds to either the quotient $1^+$ or the submodule $1^+$, but
must be lifted to the real (or at least integral) 
representation in order to implement actual mass values. In the finite representation, all we have
is a finite number of discrete place-holders for the mass.
For example, we might consider the quotient representation $1^-.1^+$ to describe a fixed charge or hypercharge, in $1^-$, 
with mass glued on top, in $1^+$. There are then three distinct vectors with a `positive' mass, namely $(0,1)$, $(1,1)$ and $(-1,1)$,
which might be used to model the three generations of any particular type of fermion.

It is also possible, though I think less likely, 
that the three generations might be better modelled by the monomial representation, that is by permutations
with a sign-change on the odd permutations. In this case, much the same is true, but with the roles of $1^+$ and $1^-$
interchanged, which may also force some changes to the interpretation.
Or, indeed, we may require the regular representation (permutations on $6$ letters) for an understanding of the interplay between mass and charge.

\subsection{Three colours of charge}
\label{3colours}
Indeed, this other PIM, $P_1^-$, which has structure $1^-.1^+.1^-$, is similar, but $1^-$ is a pseudoscalar rather than a scalar, so has a sign
attached to it. Hence it more naturally models a charge rather than a mass. The PIM is then equivalent to the representation of
$Sym(3)$ acting by permutations followed by an overall sign-change for the odd permutations. As a real representation, this
splits as $1^-+2^0$, so that, macroscopically, one can model a scalar charge or hypercharge, 
with an independent 
$2$-dimensional space of `colours'.
In the finite case, this separation of colour from charge is not possible, and one needs to consider the full $3$-dimensional
space. 

This is presumably the reason why the colours in the standard model form a $3$-dimensional space acted on by $SU(3)$,
even though the property of colour confinement suggests that a $2$-dimensional space ought to be sufficient. The finite
model would then explain why a $2$-dimensional space is not sufficient.
For a fixed generation in $1^+$, then, we have three distinct charges in $1^+.1^-$, which we might take to be $0,1,-1$ or
$-1,2/3,-1/3$ or $2/3,-1/3,-1/3$ according to context. In this way we obtain a discrete concept of `colour' inside $P_1^-$, 
acted on by $Sym(3)$, to go with the discrete concept of `generation' inside $P_1^+$.

The finite group acting on $P_1^+$ contains elements of determinant $-1$, so, when lifted to the real 
representation, cannot be embedded in a simple Lie group,
such as $SO(3)$. But it can be embedded in the unitary group $U(3)$, 
which may be why the latter group is used in the standard model.
On the other hand, the finite group acting on $P_1^-$ does embed in $SO(3)$, but at the same time the group
$GL(2,\mathbb R)$ is required for gauging the mass parameters, while the standard model combines the two
by working with a complex version of $U(2)$. 

Of course, neither $U(3)$ nor $U(2)$ can actually act on the
discrete structures $P_1^+$ or $P_1^-$, but they act on some continuous structures that mimic them. 
Indeed, these two Lie groups act on two completely different continuous structures, and therefore
they commute with each other. But in the finite model, 
the two Lie groups are replaced by a single finite group, acting simultaneously on two different discrete structures.
This creates a `mixing' between the Lie groups that is
somewhat complicated to describe, and is a significant puzzle for the standard model to explain.

\section{Representations over integers}
\label{Zreps}
\subsection{Integral representations of $Sym(2)$}
\label{S2Zreps}
The link between the modular representation theory in this paper and the ordinary representation theory in \cite{octahedral}
is provided by the integral representation theory, that is representation theory over the integer ring $\mathbb Z$. In general,
integral representation theory \cite{Reiner} is very hard, and is full of unsolved problems even for quite small groups. But the groups
$Sym(2)$ and $Sym(3)$ so far discussed have finite representation type, so it is possible to provide 
reasonably comprehensive information.

There are two irreducible representations $1^+$ and $1^-$ 
of $Sym(2)$ over $\mathbb Z$, both $1$-dimensional, where the generator $d$
of the group acts as $1$ or $-1$ respectively. There is also an indecomposable representation $2_i$, say, on $\mathbb Z^2$, in which
$d$ acts by swapping the two coordinates. In this case, $d$ has eigenvectors $(1,1)$ and $(-1,1)$ with eigenvalues $1$
and $-1$ respectively, but if we try to diagonalise $d$ with integer matrices, we have to divide by 
the determinant of the matrix of eigenvectors, which is (divisible by) $2$, so that this is impossible in $\mathbb Z$.

In other words, $d$ cannot be diagonalised, so that this $2$-dimensional representation is not $\mathbb Z$-equivalent to the direct sum
of $1^+$ and $1^-$, although they are $\mathbb R$-equivalent. 
These three representations, $1^+$, $1^-$ and $2_i$, are in fact the only indecomposable representations, so that any representation of $Sym(2)$
over $\mathbb Z$ is a direct sum of copies of them. The fact that $2_i$ is not equivalent to $1^++1^-$ may possibly be relevant to
the structure of the doublets of the weak interaction. While it is true that in the context of real or complex representation theory it makes no
practical 
difference which of the two matrices
\begin{eqnarray}
\begin{pmatrix}1&0\cr0&-1\end{pmatrix},&& \begin{pmatrix}0&1\cr1&0\end{pmatrix}
\end{eqnarray}
is used for $\gamma_5$, in the context of integral representations it may be an important distinction.

\subsection{Integral representations of $Sym(3)$}
\label{S3Zreps}
Now let us consider the group $Sym(3)$. This group has all the integral representations of its quotient $Sym(2)$, plus some
faithful representations. 
The natural permutation representation of $Sym(3)$ 
on $\mathbb Z^3$ contains a $1$-dimensional submodule spanned by $(1,1,1)$, and a $2$-dimensional
submodule spanned by $(1,-1,0)$ and $(0,1,-1)$, but the sum of these submodules only contains those $(x,y,z)$ for which $x+y+z$ is divisible by $3$.
The action of $Sym(3)$ on the $2$-dimensional submodule $2_s$ can be written as
\begin{eqnarray}
w\mapsto \begin{pmatrix}0&1\cr -1&-1\end{pmatrix}, &&
d\mapsto \begin{pmatrix}1&1\cr 0&-1\end{pmatrix},
\end{eqnarray}
while the action on the quotient $2_q$ by the $1$-dimensional submodule can be written as
\begin{eqnarray}
w\mapsto \begin{pmatrix}0&1\cr -1&-1\end{pmatrix},&&
d\mapsto \begin{pmatrix}1&0\cr-1&-1\end{pmatrix}.
\end{eqnarray}
If these two irreducible representations were $\mathbb Z$-equivalent, then the equivalence would be given by conjugation by
an integer matrix with
determinant $\pm 1$ that commutes with $w$, but the only such matrices are $\pm w$ and $\pm w^2$, so this is impossible.
Indeed, the reduction modulo $3$ tells us that one of these representations has the $3$-modular structure $1^+.1^-$, and the other
has the structure $1^-.1^+$.

A complete classification of the integral representations of $Sym(3)$ is problematic, since the Krull--Schmidt Theorem does not hold
in this context, which means that in the general case decomposition as a direct sum of indecomposables is not uniquely determined.  
But it is at least possible to classify the indecomposable representations, since all the Sylow subgroups are cyclic.
There are two indecomposable representations on $\mathbb Z^3$, corresponding to the $3$-modular indecomposables
$1^+.1^-.1^+$ and $1^-.1^+.1^-$, namely the permutation representation $3_p$ and the monomial representation $3_m$ respectively.
For reference, here is a complete list \cite{MPLee} of the $\mathbb Z$-indecomposables, with the associated real representations and 
the reductions modulo $2$ and $3$.

\begin{eqnarray}
&&\begin{array}{c|ccc}
\mathbb Z &\mathbb R & \bmod 2 & \bmod 3\cr\hline
1^+ & 1^+ & 1 & 1^+\cr
1^- & 1^- & 1 & 1^-\cr
2_i & 1^++1^- & 1.1 & 1^++1^-\cr
2_s & 2^0 & 2 & 1^+.1^-\cr
2_q & 2^0 & 2 & 1^-.1^+\cr
3_p & 1^++2^0& 1+2 & 1^+.1^-.1^+\cr
3_m& 1^-+2^0 & 1+2 & 1^-.1^+.1^-\cr
4_s& 1^++1^-+2^0& 1.1+2 & 1^-+1^+.1^-.1^+\cr
4_q& 1^++1^-+2^0& 1.1+2 & 1^++1^-.1^+.1^-\cr
6 & 1^++1^-+2^0+2^0 & 1.1+2+2 & 1^+.1^-.1^+ + 1^-.1^+.1^-\cr
\hline
\end{array}
\end{eqnarray}

\subsection{Possible interpretations}
\label{interpretations}
This is only a toy model, and there may be various physical interpretations of the different representations.
In \cite{octahedral} I suggested an interpretation in terms of three generations of leptons, but 
an interpretation in terms of all the first generation fermions might also be possible. 
It would seem reasonable to take
the odd-dimensional representations
to represent fermions, and the even-dimensional representations for bosons.

Then it would make sense to put the electron in $1^-$ and the neutrino in $1^+$, the proton in $3_p$
and the neutron in $3_m$, so that the dimension counts the number of constituent elementary particles. 
(An alternative, that does not conform to the even/odd boson/fermion principle, is to put the electron in $2^0$,
so that it has both left-handed and right-handed pieces, and then perhaps put the neutrino in $1^-$.) 
Then the question becomes, how to analyse $3_p$ and $3_m$ in terms of quarks.
The essential point is that if we just take the real constituents $1^\pm$ and $2^0$, and make them integral, then
we do not see the full structure. In the case of the permutation representation $3_p$, the constituent $1^+$ is spanned by
$(1,1,1)$, and the constituent $2^0$ by $(1,-1,0)$ and $(0,1,-1)$. To get the coordinate vectors as linear combinations
of these basis vectors, we have: 
\begin{eqnarray}
(1,0,0)&=&\textstyle\frac13(1,1,1)+\frac23(1,-1,0)+\frac13(0,1,-1),\cr
(0,1,0)&=&\textstyle\frac13(1,1,1)-\frac13(1,-1,0)+\frac13(0,1,-1),\cr
(0,0,1)&=&\textstyle\frac13(1,1,1)-\frac13(1,-1,0)-\frac23(0,1,-1).
\end{eqnarray}
In other words, it appears as though the fractional charges on the quarks arise from the fact that the direct sum structure of the
real approximation does not apply to the underlying discrete structure of the elementary particles. If we want to model the neutron in this
representation, then we can take the coefficient of $(1,-1,0)$ to be the charge of each quark. If we want to model the proton
instead, we can take the sum of the other two coefficients to be the charge.

An interpretation along these lines would seem to imply an interpretation of the $2$-dimensional representations  
as pseudoscalar mesons, hence pions.
If so, then presumably $2_i$ is the neutral pion, and $2_s$ and $2_q$ the charged pions in some order. In this way, the direct sum decomposition
of $2_i$ modulo $3$ might correspond to the standard model superposition of $u\bar u$ and $d\bar d$. On the other hand, this toy model
may not be able to distinguish particles from anti-particles, in which case we would need to interpret $2_s$ and $2_q$ as charged
and neutral pions, in some order, and perhaps leave $2_i$ for the photon, as a massless (unglued) combination of left-handed
and right-handed spins.

More speculatively, perhaps the representations $4_s$ and $4_q$ might correspond to the $W^\pm$ bosons, and $6$ to the $Z$ boson.
Or, if particle/antiparticle pairs are not distinguished, $4_s$ and $4_q$ could be allocated to both $W$ and $Z$ bosons, leaving $6$
perhaps for the Higgs boson.
Such an allocation might allow us to model the weak force for one generation in this toy model. But it also suggests a possible analysis of
the $W$ (and $Z$) bosons as containing four quarks, that is not part of the standard model. This idea is not necessarily in conflict with the standard model,
since four-quark entities do not play an important role in the latter.

\subsection{Glue}
\label{glue}
The main reason for studying the integral representation theory, however, is to study the gluing of representations.
This mathematical gluing comes in two forms: modulo $2$ and modulo $3$. The gluing modulo $3$ already appears to
correspond to the physical `gluing' of quarks into mesons and baryons. We have two dimensions of meson glue in $2_s$ and $2_q$,
corresponding to the two colourless gluons of the standard model, and we have three dimensions of baryon glue in each of
$3_p$ and $3_m$, corresponding to the six coloured gluons in the standard model. 

But the finite model reveals 
dependencies between the various gluons, so that all six ways of gluing $1^+$ to $1^-$ are mathematically equivalent,
which reduces the theoretical number of `gluons' from $8$ to $3$. This does not, of course, necessarily mean that the number of physical gluons
is $3$. It just means that the mathematics has deconstructed the gluons into the constituent colours, or some other
fundamental triplet.

The gluing modulo $2$ is quite different, and among the indecomposables of dimension $3$ or less only occurs in $2_i$.
Where the gluing modulo $3$ corresponds to $3$-dimensional indecomposables, corresponding to an action of the
standard model $SU(3)$, here we should expect an action of a standard model $SU(2)$, and therefore a
correspondence with the weak interaction. The gluing we obtain in this way seems to be a gluing between
left-handed and right-handed spinors, although, of course, in this toy model we have no spinors yet.
If we include the $4$-dimensional indecomposables $4_s$ and $4_q$, then we have three copies of this
gluing, to match up with the three dimensions of $SU(2)$. Or we could include all the indecomposables, to get
four copies of this glue altogether. Therefore we have the right number of copies to match up with the
standard model of electro-weak interactions. Whether such a matching can actually be achieved is, of course,
an entirely different question.

\section{A four-point model}
\label{4points}
\subsection{$3$-modular representations of $Sym(4)$}
\label{S4mod3}
Before we get to the spinors, we must look at the representation theory of $Sym(4)$, both over $\mathbb F_2$ and over $\mathbb F_3$. 
I will begin with the latter, since the fact that $Sym(4)$ has a normal subgroup of order
$2^2$ introduces extra complications into the $2$-modular representation theory.
The theory here is very similar to the $2$-modular representation theory of $Sym(3)$, just a bit bigger. 
In that case, all the difficulties had already been encountered in $Sym(2)$, and the extension was straightforward.
Similarly here, all the difficulties have already been encountered in $Sym(3)$.

In addition to the representations already discussed, there is a $3$-dimensional
representation coming from the rotation symmetries of a cube, and another which differs from it only in the sign of $d$.
These representations $3^\pm$ can be expressed by matrix generators
\begin{eqnarray}
x=\begin{pmatrix}1&0&0\cr0&-1&0\cr0&0&-1\end{pmatrix},&
w=\begin{pmatrix}0&1&0\cr0&0&1\cr1&0&0\end{pmatrix},&
d=\pm\begin{pmatrix}1&0&0\cr0&0&1\cr0&1&0\end{pmatrix}.
\end{eqnarray}
They 
remain irreducible when the matrix entries are regarded as elements of $\mathbb F_3$ rather than
real numbers. Hence they generate two copies of the $3\times 3$ matrix algebra as summands of the group algebra.
The remaining $6$ dimensions form the group algebra of $Sym(3)$, which I have already described in detail in Section~\ref{S3mod3}.

In other words, the algebra splits up into the $Sym(3)$ algebra already discussed in the context of the weak and strong interactions,
and two $3\times 3$ matrix algebras. 
Moreover, these matrix algebras lift to the real numbers, so that they describe quantisation of something that can be
interpreted, and therefore measured, on a macroscopic scale. For example, this could include momentum, or angular momentum, or both.
Or it could include three generations of mass, for one or two different types of particles, with different values of charge, say.

The other part of the algebra, however, describes something that only exists on the quantum scale, and cannot be lifted
to characteristic $0$. Hence this may include colours of quarks, for example, with colour confinement imposed by the fact that
$(w-1)^3=0$. All of this points to $w$ describing the fundamental quantum structure of the strong interaction, and $d$
describing the quantum structure of the weak interaction. Then the fact that $w$ and $d$ do not commute creates an
intrinsic `mixing' of the weak and strong forces, that was first introduced into the standard model with the Cabibbo angle \cite{Cabibbo}
in the early 1960s. 

The Cabibbo angle itself cannot of course be implemented in the finite representation, but requires
the real representation. It is therefore a function of how we do the reduction modulo $3$, and in particular, how we
choose a basis of particles for the two $3$-spaces. 
It is for this purpose of relating real representations to modular representations
that we need to study the integral representation theory. However, the integral representation theory of $Sym(4)$
is much harder than that of $Sym(3)$, and I cannot give a comprehensive account here.

The only other aspect of the $3$-modular representation theory of $Sym(4)$ that needs to be mentioned at this point is the
structure of the tensor products with $3^+$ and $3^-$. Since these representations are projective, all the tensor products are
projective, so the structure is completely determined by the characters. Hence all these tensor products decompose into
directs sums of copies of $3^+$ and $3^-$. 

One might interpret this property as saying that,
when we consider any measurement of a particle in the $Sym(3)$ algebra,
with any quantum property that can be lifted to real numbers, all the glue disappears from the experimental measurement. 
If so, then again we have an experimentally verified phenomenon appearing to arise naturally from the
mathematical structure of the modular representation theory.

\subsection{The Klein fours-group}
\label{V4mod2}
Before we look at the $2$-modular representations of $Sym(4)$, we need to look at the representations of its normal
subgroup, the Klein fours-group (Vierergruppe) $V_4$, consisting of elements $1,x,y,z$ that square to $1$ and satisfy
\begin{eqnarray}
xy=yx=z, \cr xz=zx=y,\cr yz=zy=x.
\end{eqnarray}

The group algebra over $\mathbb F_2$ has one $1$-dimensional ideal, three $2$-dimensional ideals and one $3$-dimensional ideal.
The $2$-dimensional ideals are generated by $1+x$, $1+y$ and $1+z$ respecitively. Any two of them intersect in the
$1$-dimensional ideal generated by $1+x+y+z$, and any two generate the $3$-dimensional ideal, which cannot be generated by
a single element.

In particular, there are three distinct $2$-dimensional indecomposable representations, distinguished by which of the three
elements
$x,y,z$ acts trivially. These three representations can be considered to be made with three different `colours' of glue,
labelled with the letters $x,y,z$. Since $xyz=1$, there is a sense in which the three colours combine to a colourless
combination (`colour confinement'). But any particular piece of glue really consists of two of $x,y,z$ acting identically.
This may therefore give us a toy model of gluons, but without the colour/anti-colour distinction, for which we need to lift to
the quaternion group. 

Now we can choose a composition series for the group algebra, as a representation of $V_4$, in three different ways,
containing the ideal generated by either $1+x$ or $1+y$ or $1+z$. Suppose we choose $1+x$. Then the bottom piece of glue
is of type $yz$, and the middle piece of glue is of a different type, but 
there is no real distinction between $xy$ and $xz$ here. 
Similarly for the top piece of glue, which loses all semblance of colour. Suppose we 
try to interpret this group algebra as a stable unit
of four elementary fermions, such as three quarks and an electron forming a hydrogen atom, 
perhaps with the electron at the top, and the three quarks on the bottom. 
Then the proton is well-defined, as the $3$-dimensional ideal, but the $2$-dimensional ideals (gluons) inside it are not.
There is therefore a choice for how to describe the internal structure of the proton, in terms of the `colours' $x,y,z$.

Now suppose we want to interpret the electro-weak interactions, by distinguishing the two up quarks from the down quark,
thereby choosing a particular $2$-dimensional ideal, say $1+x$, and the corresponding pair of colours $y,z$ on the up quarks.
While $x$ denotes a colour, which is not observable, $1+x$ attaches energy to the colour, which then becomes an
observable mass. We then have a $2$-dimensional representation consisting of two up quarks glued together,
and a $2$-dimensional quotient in which we see an electron and a down quark, each with a mass taken from the
total energy of the system.

Already we see a surprisingly large amount of the known structure of elementary particles encoded in the $2$-modular
representation theory of the Klein fours-group. In order to measure or calculate actual mass values, of course, we need to lift to
the integral representations. But there is still more information in the $2$-modular representations. In fact there are infinitely many
indecomposables, which potentially allows us to construct infinitely many `particles'. There is only one irreducible representation,
that is the trivial $1$-dimensional representation, but two irreducibles can be glued together in $3$ ways. It turns out that the gluing
can be arranged in such a way that there are indecomposables of arbitrarily large dimension \cite{Greenring}.

It may also be worth pointing out an analogy between the $16$ elements of the group algebra $\mathbb F_2V_4$ and the
$16$ dimensions of the Dirac algebra.
\begin{eqnarray}
\begin{array}{ccc}
0 & 1\cr
1,x,y,z & \gamma_0,\gamma_1,\gamma_2,\gamma_3\cr
1+x, \ldots, y+z,\ldots &\gamma_0\gamma_1,\ldots, \gamma_2\gamma_3,\ldots\cr
x+y+z,1+y+z,\ldots 
&\gamma_1\gamma_2\gamma_3, \gamma_0\gamma_2\gamma_3,\ldots\cr
1+x+y+z & \gamma_5
\end{array}
\end{eqnarray}
The correspondence is not exact, but it is suggestive. First we seem to need to lift to the quaternion group,
to obtain some of the necessary anti-commutation rules.
Then, since addition in the group algebra corresponds to multiplication in the Dirac algebra, we must
exponentiate the group algebra in some way. 
A further mixing between $1$ and $1+x+y+z$ would also seem to be required.

\subsection{Triplet symmetries}
\label{triplets}
Lifting from $Sym(3)$ to $Sym(4)$ splits the basic triplet symmetry $w$ into four distinct triplet symmetries,
$w$, $xw$, $yw$ and $zw$.
This gives us some scope for modelling the various different triplet symmetries in the standard model.
These include colour symmetry, generation symmetry (probably split further into lepton generations and quark generations),
and the up/down/strange symmetry that underlies the original Eightfold Way.

It is not clear whether these are all distinct in the finite $Sym(4)$ model. 
All of them act by conjugation to cycle $x,y,z$, but they can be distinguished by multiplication, or equivalently by
the action as permutations on $4$ letters:
\begin{eqnarray}
w&\mapsto &(X,Y,Z),\cr
xw&\mapsto &(W,Y,X),\cr
yw&\mapsto &(W,Z,Y),\cr
zw&\mapsto &(W,X,Z).
\end{eqnarray}
Thus $w,xw,yw,zw$ fix $W,Z,X,Y$ respectively. The letters $W,X,Y,Z$ are arbitrary, so we might as well
use $W$ for energy/mass, so that $w$ would be a colour symmetry, that does not affect the mass.
Then we could use the other three triplet symmetries
for lepton generations, quark generations and up/down/strange symmetries, for example.

There are, of course, other possible interpretations, but in all cases the triplet symmetries should tell us something about mass,
not of individual particles, but of triplets of particles. That is, the sum of the masses of three carefully chosen particles should have
a more fundamental character than the individual masses. The type for such particle triplets is set by the Coleman--Glashow relation
\cite{ColemanGlashow}, that uses a rotation of up to down to strange quarks to map the approximate equality of proton and neutron
masses to approximate equalities between $\Xi^0$ and $\Xi^-$ baryons, and between $\Sigma^+$ and $\Sigma^-$, and by summing over
the triplets obtains an `exact' \cite{CGexact} mass equality
\begin{eqnarray}
p+\Sigma^-+\Xi^0&=& n+\Sigma^++\Xi^-.
\end{eqnarray}

This equation suggests that a similar equation might be obtained from the approximate equality of the neutron mass to the sum of
the electron and proton masses, using the lepton generation symmetry instead of the up/down/strange symmetry. Assuming that the
proton is not affected by the lepton generation symmetry, the total mass on one side of the equation comes from three generations
of electrons and three protons. Clearly this is nowhere near the mass of three neutrons. But it turns out to be equal to the
mass of \emph{five} neutrons, to well within current experimental uncertainty \cite{perspective}:
\begin{eqnarray}
\label{emutau}
e+\mu+\tau + 3p &=& 5n.
\end{eqnarray}
The model does not explain where this number $5$ comes from, but perhaps it is a clue to a \emph{finite} unification of
weak doublets with strong triplets, analogous to, but quite different from, the Georgi--Glashow model \cite{GG}?
Or is it just a meaningless coincidence?

I shall return to this question in Section~\ref{mass}, where I hope to show that
a detailed look at the spinors provides some further justification for this equation.
Essentially what the model suggests is to count $15$ fundamental units on both sides of the equation,
by splitting the leptons into left-handed and right-handed units, and spitting the baryons into
red, green and blue units.
Both these triplet equations are consistent with experiment, to a current accuracy of around one part in $10^4$.
It is possible, therefore, that they are exact, and follow from the discrete structure of the groups and their representations.
\subsection{Further symmetries}
\label{doublets}
There are also a number of approximate mass equations, which cannot therefore arise from this discrete structure,
so that we must look elsewhere for an explanation. These include the near equality of proton and neutron masses,
with a relative difference of around $14\times 10^{-4}$, and of the $\Xi^0$ and $\Xi^-$ baryons ($50\times 10^{-4}$),
in both cases exhibiting a doublet symmetry.

Another interesting case comes from the quadruplet symmetry of the $\Lambda$ and $\Sigma$ baryons.
All these baryons have a strangeness of $-1$, and the sum of their masses is very close to five neutron masses.
But the (tempting) hypothesis that these masses are exactly equal is not consistent with experiment. The relative difference is around
$5\times 10^{-4}$, of the same order of magnitude as the proton/neutron mass difference, so perhaps (if it isn't just a random coincidence) 
arises from
the same fundamental cause, whatever that may be.

The ultimate mathematical reason for the triplet mass equations may lie in the $3$-modular representation theory, in such a way that
once we have chosen the particular triplet symmetry we are interested in, the symmetry group is reduced to $Sym(3)$. If so, then
doublet and quadruplet symmetries may similarly reduce the symmetry group to $V_4$ or $D_8$, which suggests a role for the
$2$-modular representation theory in this case. The physics then suggests that we are in the realm of the weak interaction,
in which energy is a scalar but mass is not, rather than 
the strong interaction, in which mass is a scalar and energy, therefore, is not. In other words, in this situation we should
not expect exact equality of masses. Nevertheless, the very close matching of unequal masses in these situations still
begs for some kind of explanation.

Another possible interpretation of the four basic triplet symmetries is as generation symmetries for the four different types of
fundamental fermions: neutrinos, electrons, and up and down quarks. This would give the neutrinos a theoretical mass of zero,
but this is not inconsistent with a small \emph{effective} mass for neutrinos in a gravitational field: since mass is defined in the
integral group ring, or even in the real group algebra, it can be taken to be defined by a vector very close to, but not exactly equal to, $W$.
At the discrete level, however, we would have to take the neutrino mass to be `intrinsically' zero, and identify the neutrino
generation symmetry with the colour symmetry of quarks.
 
 Of course, this should not be interpreted as saying that neutrino generation `is' colour, only that the symmetry groups of the two
 concepts are equal. But this does mean that the symmetry group of the gluons would then be essentially the same as the symmetry group of the
 neutrino/anti-neutrino pairs, modulo scalars. It is therefore not entirely inconceivable that a theory of the mixing of the weak and strong forces might
 be built on some kind of equivalence between a gluon and a neutrino/anti-neutrino pair. 
 
 At its most radical, this idea could even
 lead to the gluons becoming redundant as mediators for the strong force, if they could be replaced in the theory by pairs
 consisting of a neutrino and an anti-neutrino. This need not change the theory of the strong interaction in any way, and would only provide
 an extra level of theoretical underpinning for the theory that already exists. But it could also provide a significant conceptual
 simplification, if it turned out to be possible to explain the strong force in terms of virtual neutrinos and anti-neutrinos,
 as an alternative to gluons. Further remarks on this speculative suggestion are made in Section~\ref{blackholes}.
 
 \section{Doubling the four points}
 \label{8points}
\subsection{The group $GL(2,3)$}
\label{GL23}
Let us now move on to considering the spinors, and how they behave modulo $3$.
The group that I studied in \cite{octahedral} happens to be isomorphic to the group $GL(2,3)$ of all invertible $2\times 2$
matrices over the field $\mathbb F_3$ 
of order $3$. It would therefore seem to be a suitable finite group to use 
to quantise some or all of the various $2$-dimensional groups
in the standard model, such as the spin group $SL(2,\mathbb C)$ and/or the electroweak gauge group $U(2)$. The $2$-dimensional
representations can then be used for quantised spinors, and the $1$-dimensional and $3$-dimensional representations
can be used to quantise bosons of various kinds.

For compatibility of notation, I take generators $i,w,d$ defined by
\begin{eqnarray}
i:=\begin{pmatrix}0&1\cr -1&0\end{pmatrix},&
w:= \begin{pmatrix}1&1\cr0&1\end{pmatrix},&
d:=\begin{pmatrix}1&0\cr 0&-1\end{pmatrix}.
\end{eqnarray}
and define
\begin{eqnarray}
j:=w^{-1}iw, && k:=wiw^{-1}.
\end{eqnarray}
It is then easy to verify the fundamental relations
\begin{eqnarray}
ij=k, & jk=i, & ki=j,\cr
i^2=-1,& j^2=-1,& k^2=-1,\cr
iw=wj, & jw=wk, & kw=wi,\cr
id=-di, & jd=-dk, & kd=-dj,\cr
wdw=d, &w^3=1, & d^2=1.
\end{eqnarray}

There are altogether six irreducible representations over $\mathbb F_3$, two each in dimensions $1$, $2$ and $3$.
Writing down only the traces of the matrices loses vital information about the multiplicities of eigenvalues, so that
one uses instead
the so-called Brauer character table \cite{ABC}, that contains complex numbers sufficient to 
recover all the eigenvalue information. The elements $w$ and $w^{-1}$ have all eigenvalues equal to $1$, so can be omitted from the
table:
\begin{eqnarray}
\begin{array}{c|cccccc}
&1&-1&i&d&jd&-jd\cr\hline
1^+ & 1 & 1 & 1 & 1&1&1\cr
1^- & 1 & 1 & 1 & -1&-1&-1\cr
2^+ & 2 & -2 & 0 & 0 & \sqrt{-2} & -\sqrt{-2}\cr
2^- & 2 & -2 & 0 & 0 & -\sqrt{-2} & \sqrt{-2}\cr
3^+ & 3 & 3 & -1 & 1 & -1 & -1\cr
3^- & 3 & 3 & -1 & -1 & 1 & 1\cr
\hline
\end{array}
\end{eqnarray}
Notice that in $\mathbb F_3$, we replace $-2$ by $+1$, so that $\sqrt{-2}$ corresponds to a trace of $\sqrt1=\pm 1$, but we are allowed to
choose our sign convention here.

Representations corresponding to these irreducible Brauer characters can be taken as follows:
\begin{eqnarray}
&&\begin{array}{c|ccc}
&i
&w & d\cr\hline
1^\pm &(1)
& (1) &\pm(1)\cr
2^\pm & 
\begin{pmatrix}0&1\cr -1&0\end{pmatrix}&
 \begin{pmatrix}1&1\cr0&1\end{pmatrix}&
\pm\begin{pmatrix}1&0\cr 0&-1\end{pmatrix}\cr
3^\pm &\begin{pmatrix}1&0&0\cr0&-1&0\cr0&0&-1\end{pmatrix}
&\begin{pmatrix}0&1&0\cr 0&0&1\cr 1&0&0\end{pmatrix} & \pm\begin{pmatrix}1&0&0\cr0&0&1\cr 0&1&0\end{pmatrix}\cr
\end{array}
\end{eqnarray}
As usual, even-dimensional representations are fermionic, and odd-dimensional representations are bosonic.
The bosonic representations are representations of $Sym(4)$, and have already been discussed in some detail (Section~\ref{S4mod3}).
The irreducibles split into a scalar $1^+$, a pseudoscalar $1^-$, a vector $3^-$ and an axial vector $3^+$, but it is the interplay between
the scalar (quantised mass) and the pseudoscalar (quantised charge) that is the interesting part of the representation theory.

The fermionic irreducible representations are two types of spinors, which play the roles of the left-handed and right-handed Weyl spinors
in the standard model. However, as I have already discussed in the context of 
the toy model,
the relationship between the concepts of `left-handed' and `right-handed' in the finite 
and the continuous cases is not entirely straightforward,
and we need to be careful with the interpretations here.
Notice in particular that, as continuous representations, they are complex Weyl spinors, not real, but as finite representations,
they do not require any extension of the ground field, since there is no requirement to represent $\sqrt{-1}$ or $\sqrt{2}$, but only
$\sqrt{-2}$.  
Thus the difference between left-handed and right-handed spinors is not just a conventional difference
between a representation and its complex conjugate, but a genuine difference between $+1$ and $-1$.

All six of the irreducible $3$-modular representations arise from reducing real representations modulo $3$, or, in the case of the spinor
representations, reducing complex representations (Weyl spinors) modulo $\sqrt{-3}$. There are two other irreducible real representations,
of dimensions $2$ and $4$, whose characters modulo $3$ are $1^++1^-$ and $2^++2^-$ respectively. The reductions modulo $3$ can be
glued either way up, depending on which integral representation is chosen for the reduction. 

\subsection{Quantisation of spacetime}
\label{quantumspacetime}
The former case has already been
discussed in the context of the representations of $Sym(3)$, while the latter was used in \cite{octahedral} for a quantisation
of spacetime, with an explicit copy given by representing $i,w,d$ with the matrices
\begin{eqnarray}
\begin{pmatrix}0&1&0&0\cr -1&0&0&0\cr 0&0&0&-1\cr 0&0&1&0\end{pmatrix},&
\begin{pmatrix}1&0&0&0\cr 0&0&1&0\cr 0&0&0&1\cr 0&1&0&0\end{pmatrix},&
\begin{pmatrix}1&0&0&0\cr 0&-1&0&0\cr 0&0&0&-1\cr 0&0&-1&0\end{pmatrix}.
\end{eqnarray}
If we label the coordinates $W,X,Y,Z$ respectively, then these matrices are obtained from the permutations
$(W,X)(Y,Z)$, $(X,Y,Z)$ and $(Y,Z)$ by adding in some signs. 

If we now interpret the matrix entries modulo $3$,
that is, in the field $\mathbb F_3$, then the representation is no longer irreducible.
Indeed, there is a submodule of dimension $2$, whose non-zero vectors are
\begin{eqnarray}
\pm(X+Y+Z),&\pm(W-X+Y),&\pm(W-Y+Z),\quad\pm(W-Z+X).
\end{eqnarray}
We have
\begin{eqnarray}
(W-X+Y)-(W-Y+Z)&=& -X +2Y-Z,
\end{eqnarray}
so that modulo $3$ these vectors form a subspace, while over $\mathbb R$ they span the whole $4$-space.
Quantum spacetime in this model, therefore, has a very different structure from macroscopic spacetime.
If we take $W-Y+Z$ and $X+Y+Z$ as basis vectors for this $2$-space over $\mathbb F_3$, 
then the matrices for $i,w,d$ are those given earlier for the
representation $2^+$.  There is then a quotient module equivalent to $2^-$. But there is no submodule $2^-$, and no quotient
module $2^+$. 

With these conventions, $2^+$ appears to denote a `left-handed' spinor, and $2^-$ denotes a `right-handed' spinor, which cannot exist
on its own but only glued on top of a left-handed spinor. We can however change our conventions to reverse these signs if we wish, essentially by
changing the sign of $W$, that is of the `time'. Interpreting this physically seems to suggest that, at least as far as the strong force is concerned,
quantised spacetime is essentially the same thing as spin. Whether this is still true in the context of electro-weak interactions, that is
on reduction modulo $2$, is a completely different question.

To look at the weak force in the $3$-modular context, we need to look at the element $id$, as the appropriate analogue of $\pm\gamma_5$
in the Dirac algebra. Now $-id$ acts by swapping $W$ and $X$, and negating $Y$, so swaps the two basis vectors of $2^+$. In the real
representation, it fixes $W+X$ and $Z$, and negates $W-X$ and $Y$. 
Thus the weak force is intimately connected to the
macroscopic directions $X,Y,Z$ that describe the embedding of quantised spin in macroscopic spacetime. In particular, the weak force is
chiral in the sense that if we negate the $X$ direction then we must also swap $Y$ and $Z$.

It is worth pointing out that this chirality is a direct consequence of the quantisation of spacetime in terms of the $3$-modular representation
theory of the $3$-modular equivalent of $SU(2)$. But there is nothing in the mathematics to say whether this chirality is left-handed or
right-handed in the physical world. It is just a mathematical convention, so that the physical chirality must come from somewhere else.

\subsection{The $3$-modular group algebra of $GL(2,3)$}
\label{groupalgebra}
The finite version of the group algebra, that is $\mathbb F_3.GL(2,3)$, 
does not have the simple structure as a direct sum of matrix algebras that the real and complex
group algebras have. It does split as a direct sum of 
a bosonic and a fermionic algebra, each of dimension $24$. The bosonic algebra is the
group algebra of $Sym(4)$ discussed in Section~\ref{S4mod3} above, and splits as a direct sum of two $3\times 3$ matrix algebras and the
$6$-dimensional group algebra of $Sym(3)$. 

The fermionic algebra does not split as a direct sum of smaller algebras,
but behaves like a scaled-up version of the group algebra of $Sym(3)$,
in which the $1\times 1$ matrices (scalars and pseudoscalars) are replaced by
$2\times 2$ matrices (left-handed and right-handed spinors). But the translation
from $Sym(3)$ is far from straightforward, and the interpretation of the
Weyl spinors is fraught with difficulty.

In particular, the group algebra splits into four $3$-blocks, one built from
scalars $1^+$ and pseudoscalars $1^-$, one from vectors $3^-$, one from axial vectors $3^+$,
and one from spinors $2^+$ and $2^-$.
Every representation splits as a direct sum of pieces from each block.
But we cannot in general
split the scalars from the pseudoscalars, or the left-handed spinors from the right-handed spinors. The scalars and pseudoscalars
can be glued on top of each other to make towers of height at most $3$, and similarly for the spinors:
\begin{eqnarray}
\begin{matrix}1^+\cr1^-\cr1^+\end{matrix},\qquad
\begin{matrix}1^-\cr1^+\cr1^-\end{matrix},&&
\begin{matrix}2^+\cr2^-\cr2^+\end{matrix},\qquad
\begin{matrix}2^-\cr2^+\cr2^-\end{matrix}.
\end{eqnarray}

We therefore have three different types of spinors, of heights $1$, $2$ and $3$. Those of height $1$ are finite versions of Weyl spinors,
so may represent neutrinos and antineutrinos, although the fundamental difference between scalars and pseudoscalars
suggests that we may need a more subtle interpretation. 
Those of height $2$ are finite versions of Dirac spinors, with the left-handed and right-handed parts
glued together into an inseparable whole, so might represent electrons, if these are not
already taken care of by height $1$ spinors.  
Those of height $3$ might represent baryons, made somehow out of three $2$-dimensional quarks `glued' together
into a $6$-dimensional whole.

Whatever interpretation we choose, we must deal with the fact that the structure of this algebra differs
markedly from the corresponding real algebra, which is the direct sum of a $2\times 2$ complex matrix algebra
(representing spinors) and a $4\times 4$ real matrix algebra (used in \cite{octahedral} to 
represent the product of spacetime with its dual). In other words,
at the quantum level the model requires a fundamental mixing together of matter (represented by the spinors)
and spacetime itself. Such a mixing is not, of course, part of the standard model of particle physics. 

But it is a
fundamental part of general relativity \cite{thooft,GR1,GR2}. Therefore we should expect to have to implement such a mixing in any 
attempt to quantise gravity, 
that reduces to general relativity in the appropriate limit.
The fact that the modular representation theory provides such a mixing  
is a strong hint that this might be
a useful way to proceed. If so, then the fundamental physical mechanism for such a quantum gravity would have to correspond
mathematically to a decoupling of the spinors from spacetime, which would seem to mean converting the spinors
into neutrinos. Of course, this is highly speculative, but not completely unreasonable, given the significant mysteries that 
still surround the neutrinos.

For example, one could 
perhaps consider the primary function of the gluons to be the coupling of matter to spacetime, 
instead of the usual interpretation of coupling of matter to matter, so that on decoupling, the gluon `decays' into a neutrino and
an anti-neutrino, which then travel through the universe and interact very weakly with other matter, in order to transmit
a very weak gravitational force. Or perhaps it is better to keep the usual interpretation of the gluons, and interpret the 
individual neutrinos and anti-neutrinos
as coupling matter to spacetime.

\section{Spinors}
\label{spinors}
\subsection{The faithful $3$-block}
\label{faithfulblock}
There are two PIMs in the faithful block, that contains the spinors. These have shape $2^+.2^-.2^+$ and $2^-.2^+.2^-$ and may
represent some kind of extension of a $3$-modular version of the Dirac spinor, to include three Weyl spinors rather than just two.
The whole block has dimension $24$, and consists of two copies of each of these PIMs. 
As $3$-modular representations they can be constructed in multiple ways as tensor products of the representations that we have
already studied:
\begin{eqnarray}
\begin{matrix}2^+\cr2^-\cr 2^+\end{matrix}&=&
\begin{matrix}1^+\cr1^-\cr 1^+\end{matrix}\otimes 2^+ =
\begin{matrix}1^-\cr1^+\cr 1^-\end{matrix}\otimes 2^- = 3^-\otimes 2^+ = 3^+\otimes 2^-,
\end{eqnarray}
and the same with $2^+$ and $2^-$ interchanged. The equality of these various tensor products may represent some kind of
$3$-modular analogue of the Dirac equation, if, as previously suggested, $3^-$ is related to momentum, $1^+$ to energy/mass and
$1^-$ to charge.

Since the characters of the PIMs are complex,
there is no particularly easy way to construct them from real or integral representations. 
But the integral representations are what we need in order to relate the PIMs to physical properties that can
be measured. 
The smallest suitable integral representation
 is $12$-dimensional, and can be constructed from the monomial representation on $\pm W, \pm X, \pm Y, \pm Z$ given above, as follows.
First take the set of $12$ vectors and their negatives:
\begin{eqnarray}
\begin{array}{c|ccc}
X-Y-Z & -W-Y+Z &-W-Z+X & -W-X+Y\cr
-X-Y+Z & W+Y+Z & W+Z+X & W+X+Y\cr
-X+Y-Z & W-Y-Z & W-Z-X & W-X-Y
\end{array}
\end{eqnarray}
Then define $12$ corresponding basis vectors for a $12$-dimensional space, with the labels
\begin{eqnarray}
\begin{array}{c|ccc}
t_0 & x_0&y_0&z_0\cr
t_+ & x_+&y_+&z_+\cr
t_- & x_-&y_-&z_-
\end{array}
\end{eqnarray}
The notation is chosen so that the triplet 
$x,y,z$ is closely related to the triplet $X,Y,Z$, and the subscripts $0,+,-$ encode the signs.
I have used $t$ instead of $w$, since $w$ is already in use for something else. In fact, the symbols $t,x,y,z$
used here will turn out to have some relationship to discrete spacetime coordinates, although this relationship
is not straightforward to describe.

The actions of the group generators on this $12$-space are given by permutations with signs attached:
\begin{eqnarray}
w&=&(t_0,t_-,t_+)(x_\alpha,y_\alpha,z_\alpha)\cr
d&=&
(t_\alpha,-t_{-\alpha})(x_\alpha,x_{-\alpha})(y_\alpha,z_{-\alpha})\cr
i&=&(t_\alpha,x_\alpha,-t_\alpha,-x_\alpha)(z_\alpha,y_{\alpha-1},-z_\alpha,-y_{\alpha-1})
\end{eqnarray}
Here, arithmetic on the subscripts is carried out modulo $3$. Ignoring the subscripts we see an action of $w$ rotating the internal space
directions $x,y,z$, while $d$ splits into a T symmetry negating $t$ and a P symmetry swapping $y$ with $z$.

\subsection{Subspaces}
\label{subspaces}
Within this $12$-space there are many interesting invariant subspaces spanned by certain integral linear combinations of the basis vectors.
For example, if we define
\begin{eqnarray}
A=\begin{array}{|c|ccc|} \hline.&+&+&+\cr.&.&.&.\cr.&.&.&.\cr\hline\end{array},&&
B=\begin{array}{|c|ccc|} \hline+&.&.&.\cr.&.&.&+\cr.&.&-&.\cr\hline\end{array},\cr
C=\begin{array}{|c|ccc|} \hline.&.&.&.\cr.&+&.&.\cr+&.&.&-\cr\hline\end{array},&&
D=\begin{array}{|c|ccc|} \hline.&.&.&.\cr+&.&+&.\cr.&-&.&.\cr\hline\end{array},
\end{eqnarray}
which is just a pictorial representation of the equations
\begin{eqnarray}
A:=x_0+y_0+z_0,&& B:=t_0-y_-+z_+, \cr &&C:=t_--z_-+x_+, \cr && D:=t_+-x_-+y_+,
\end{eqnarray}
then the action on $A,B,C,D$ is the same as the action on $W,X,Y,Z$.
In other words, this subspace has a potential interpretation as macroscopic spacetime, within the larger space that also contains
$8$ more dimensions of spinors of various kinds.
On reduction modulo $3$, the $4$-space spanned by $A,B,C,D$ has an invariant subspace spanned by $A-C+D$ and $B+C+D$, on which
the group acts via the matrices given above for the representation $2^+$. It follows that the whole $4$-space modulo $3$ is indecomposable
and has the structure 
$2^+.2^-$. 

There is another $4$-space in the integral representation that reduces modulo $3$ to $2^-.2^+$. An example is given by
the spanning vectors
\begin{eqnarray}
P=\begin{array}{|c|ccc|} \hline+&.&.&.\cr+&.&.&.\cr+&.&.&.\cr\hline\end{array},&&
Q=\begin{array}{|c|ccc|} \hline.&+&.&.\cr.&+&.&.\cr.&+&.&.\cr\hline\end{array},\cr
R=\begin{array}{|c|ccc|} \hline.&.&+&.\cr.&.&+&.\cr.&.&+&.\cr\hline\end{array},&&
S=\begin{array}{|c|ccc|} \hline.&.&.&+\cr.&.&.&+\cr.&.&.&+\cr\hline\end{array},
\end{eqnarray}
or more explicitly
\begin{eqnarray}
P:=
t_0+t_++t_-,&& Q:=x_0+x_++x_-,\cr
&&R:=y_0+y_++y_-,\cr
&&S:=z_0+z_++z_-,
\end{eqnarray}
on which the group generators act as
\begin{eqnarray}
w&=&(Q,R,S)\cr
d&=&(P,-P)(R,S)\cr
i&=&(P,Q,-P,-Q)(R,-S,-R,S).
\end{eqnarray}
Here the subspace modulo $3$ is spanned by $P-R+S$ and $Q+R+S$.

As real representations, $A,B,C,D$ and $P,Q,R,S$ span two copies of the $4$-dimensional irreducible. 
If we simply add up the corresponding vectors in the space spanned by $W,X,Y,Z$, then we see that
\begin{eqnarray}
A,B,C,D &\mapsto& 3W,3X,3Y,3Z\cr
P&\mapsto & X+Y+Z\cr
Q&\mapsto & -W+Y-Z\cr
R&\mapsto & -W+Z-X\cr
S&\mapsto & -W+X-Y
\end{eqnarray}
Hence it would appear that these two copies of the same real representation must have different physical
interpretations. If we want to interpret $A,B,C,D$ as macroscopic spacetime, then we would seem to need to interpret
$P,Q,R,S$ as some type of spin. Before we try to interpret this spin in terms of the standard model left-handed and
right-handed spins, let us look at the rest of the spin representation in more detail.

\subsection{An integer version of the Weyl spinor}
\label{Zspinor}
The rest of the real $12$-space
consists of a Majorana--Weyl spinor, and maps to zero in the space spanned by $W,X,Y,Z$.
It can be obtained as the orthogonal complement of the $8$-space already exhibited.
It is therefore easy to find a basis for this real $4$-space consisting of the vectors
\begin{eqnarray}
H=\begin{array}{|c|ccc|} \hline.&+&-&.\cr-&-&+&.\cr+&.&.&.\cr\hline\end{array},&&
I=\begin{array}{|c|ccc|} \hline-&.&+&-\cr+&.&-&+\cr.&.&.&.\cr\hline\end{array},\cr
J=\begin{array}{|c|ccc|} \hline.&+&.&-\cr-&.&.&.\cr+&-&.&+\cr\hline\end{array},&&
K=\begin{array}{|c|ccc|} \hline+&.&-&+\cr.&.&.&.\cr-&.&+&-\cr\hline\end{array},
\end{eqnarray}
and to compute the action of the group as integral matrices acting on this basis.

 It is possibly slightly more revealing to change to a basis such as
$I+K,J-H,H,I$ so that the matrices exhibit certain block structures, although we should not expect too much
since this representation is irreducible.
\begin{eqnarray}
i&=&\begin{pmatrix}0&1&&\cr-1&0&&\cr0&-1&0&1\cr -1&0&-1&0\end{pmatrix},\cr
w&=&\begin{pmatrix}0&1&0&-1\cr-1&-1&-1&0\cr &&0&1\cr&&-1&-1\end{pmatrix},\cr
d&=&\begin{pmatrix}1&&&\cr0&-1&&\cr0&1&1&\cr 1&0&0&-1\end{pmatrix}
\end{eqnarray}

Now it is possible to convert this real representation into a $2$-dimensional complex representation,
by first computing the class sum on one of the classes of elements of order $8$, such as
\begin{eqnarray}
(j-k)d+(k-i)wd+(i-j)w^2d.
\end{eqnarray}
Since this class sum represents the scalar $3\sqrt{-2}$ on one of the $2$-dimensional complex representations, this allows us to
define $\sqrt{-2}$.

The calculations are not so easy to carry out correctly by hand, but 
it turns out that (up to a conventional choice of sign) we have
\begin{eqnarray}
\sqrt{-2}&:=& \begin{pmatrix} 1&1&2&0\cr-1&-1&0&2\cr -1&0&-1&-1\cr 0&-1&1&1\end{pmatrix},
\end{eqnarray}
so that we can take $H,I$ to form a complex basis with the definitions
\begin{eqnarray}
\sqrt{-2}H &=& -H-K-2I,\cr
\sqrt{-2}I&=&I-J+2H.
\end{eqnarray}
Then it is easy to write down the corresponding complex matrices for the group generators:
\begin{eqnarray}
i=\begin{pmatrix}-1&\sqrt{-2}\cr\sqrt{-2}&1\end{pmatrix},&&
w=\begin{pmatrix}0&1\cr-1&-1\end{pmatrix},\cr
&&d=\begin{pmatrix}2&1-\sqrt{-2}\cr -1-\sqrt{-2}&-2\end{pmatrix}.
\end{eqnarray}

\subsection{Extension to the Dirac spinor}
\label{DiracZspinor}
At this stage we have split the real $12$-space into 
three distinct $4$-spaces,  one of which appears to transform like spacetime, and two of which
look like spinors of some kind.
The two spinors, however, are quite different.
One of them has a natural complex structure, so can plausibly be interpreted as a Majorana--Weyl spinor,
while the other one does not. So although it seems reasonable to combine these two spinors into a Dirac spinor,
the resulting complex structure is an artificial construct, and does not respect the underlying 
discrete structure of the model.

In other words, if we make a choice of complex structure on the spinor spanned by $P,Q,R,S$ in order to
map to a standard model Weyl spinor, then we lose
some of the underlying symmetry. Since we really only have a choice of three complex structures, 
this is a plausible point at which the standard model has been forced to restrict to a single generation.
To compute these complex structures, which apply only to a subgroup of index $3$ in $GL(2,3)$, we may define
one of the following three elements to be $\pm\sqrt{-2}$:
\begin{eqnarray}
(j-k)d, &(k-i)wd,&(i-j)w^2d.
\end{eqnarray}
First we calculate
\begin{eqnarray}
jd&=&(P,S,-Q,-R,-P,-S,Q,R),\cr
-kd&=&(P,-R,Q,S,-P,R,-Q,-S)
\end{eqnarray}
from which we read off, for one of the two possible choices of sign,
\begin{eqnarray}
\sqrt{-2}P = S-R,&&\sqrt{-2}R=P-Q,\cr
\sqrt{-2}Q= S+R,&&
\sqrt{-2}S=-P-Q.
\end{eqnarray}
Notice, incidentally, that $P$ and $Q$ do not form an integral basis for the representation, since obtaining $R$ and $S$ from these
requires dividing by $\sqrt{-2}$. Either $P,R$ or $P,S$ will do instead.
In this way we obtain an integral basis such as $H,I,P,R$ for a complex $4$-space, that would appear to be a discrete
version of the Dirac spinor for a single generation of fermions. 

The other two choices of complex structure can be obtained by cycling $Q,R,S$, and give spinors that we should expect to
be applicable to the other two generations of fermions. However, we should be careful here, since the generations of quarks do not match
to the generations of leptons, and perhaps restrict our attention to lepton generations at this stage.
Indeed, the phenomenon of neutrino oscillations \cite{oscillation,neutrinos,SNO}
suggests that we should also be careful not to talk about neutrino generations
until the picture becomes clearer than it is at present. Hence all we can really say here is that the choice of complex structure
corresponds to a choice of the generation of the electron.

\section{The standard model and beyond}
\label{SMplus}
\subsection{The Dirac matrices}
\label{gammas}
In \cite{octahedral} I approximately identified the Dirac matrices $\gamma_1$, $\gamma_2$, $\gamma_3$ and $\gamma_0$ with the elements
$i$, $j$, $k$ and $d$ of the group, although these elements do not satisfy quite the same relations as the Dirac matrices. In addition
I suggested using $id$ in place of $\gamma_5$. In order to see how closely this matches, or does not match, the standard model,
we need to calculate the complex matrices corresponding to $i,j,k,d$ and $id$, on the two complex $2$-spaces spanned by $H,I$ and $P,R$.
The latter is much easier to calculate, and gives
\begin{eqnarray}
\label{DiracPR}
i=\begin{pmatrix}1&-\sqrt{-2}\cr -\sqrt{-2}&-1\end{pmatrix}, & j=\begin{pmatrix}0&1\cr-1&0\end{pmatrix},&
k=\begin{pmatrix}\sqrt{-2}&1\cr1&-\sqrt{-2}\end{pmatrix},\cr
d=\begin{pmatrix}-1&0\cr\sqrt{-2}&1\end{pmatrix},&id=\begin{pmatrix}1&-\sqrt{-2}\cr 0&-1\end{pmatrix},&
1-id=\begin{pmatrix}0&\sqrt{-2}\cr 0&2\end{pmatrix}.
\end{eqnarray}
In particular, we clearly see the projection with $1-id$ onto the imaginary and real parts of the second coordinate of the complex $2$-space.

In the integer version, however, there is a factor of $\sqrt{-2}$ which we might need to take out. In other words, the standard model
projection with $(1-\gamma_5)/2$ might correspond to $(1-id)/\sqrt{-2}$, or perhaps it is better to say that $1-\gamma_5$ corresponds to
$\sqrt{-2}(1-id)$. In any case, this factor of $\sqrt{-2}$ has some subtle effects which we will need to investigate. It may play a role
similar to that played by the imaginary factor in the standard model equation $\gamma_5=i\gamma_0\gamma_1\gamma_2\gamma_3$.

On the complex $2$-space spanned by $H$ and $I$ we have
\begin{eqnarray}
\label{DiracHI}
i=\begin{pmatrix}-1&\sqrt{-2}\cr \sqrt{-2}&1\end{pmatrix},&&
j=\begin{pmatrix}1+\sqrt{-2}&2\cr -\sqrt{-2} &-1-\sqrt{-2}\end{pmatrix},\cr
k=\begin{pmatrix}1-\sqrt{-2}&-\sqrt{-2}\cr -2 & -1+\sqrt{-2}\end{pmatrix},&&
d=\begin{pmatrix}2&1-\sqrt{-2}\cr -1-\sqrt{-2} & -2\end{pmatrix},\cr
id=\begin{pmatrix}-\sqrt{-2} & -1-\sqrt{-2}\cr -1+\sqrt{-2} & \sqrt{-2}\end{pmatrix},&&
1-id=\begin{pmatrix}1+\sqrt{-2} & 1+\sqrt{-2}\cr 1-\sqrt{-2} & 1-\sqrt{-2}\end{pmatrix}.
\end{eqnarray}
Here again we see a projection, since the two columns of $1-id$ are equal. Over the complex numbers,
we can change basis in order to diagonalise $id$ in both representations, so that everything matches up,
by matching the eigenvectors of $id$:
\begin{eqnarray}
R & \leftrightarrow &H+I\cr
S = \sqrt{-2}P + R &\leftrightarrow&\sqrt{-2}(H+I) + (I-H) = H-I-J-K.
\end{eqnarray}

 Over the
integers, however, it is impossible to match up the two halves of the spinor.
There are some very significant differences between the $H$-spinor and the $P$-spinor, which
are not so easy to explain away.
Indeed, I would suggest that it is better not to try to explain them away, but to use the extra subtleties of the
integer representations to incorporate some of the subtleties of the standard model, that cannot be incorporated in
the complex Dirac spinor.

The part of the spinor that we appear to need to identify with the standard model `left-handed' part is then
spanned by $R$ and $H+I$, if we continue to write the spinors as rows. It is, of course, not obvious that this identification
makes sense. It may well be that the proposed changes to the Dirac spinor break something important, rather than
adding something new. 
But making such a choice of $R$ and $H+I$ breaks the symmetry of the finite group, which has already been reduced from
$GL(2,3)$ to the semi-dihedral group $SD_{16}$, down to nothing but the scalar $-1$.

\subsection{Continuous structures}
\label{continuous}
It is now time to start to introduce the continuous symmetries that are used in the standard model. Up to this point, I have used
$\gamma_1$, $\gamma_2$ and $\gamma_3$ to model discrete symmetries, whereas they are used in the standard model to model
continuous symmetries. To be more precise, they are used for momentum, which is quantised in magnitude but not in direction, so that,
in principle at least, there
is one discrete variable and two continuous variables here. Or, to be even more precise, they are used together with $\gamma_0$
in order to include mass, which is also not quantised, so that altogether there is one discrete variable and three continuous variables.
Moreover, the discrete versions of $\gamma_1$, $\gamma_2$ and $\gamma_3$
do not satisfy the same relations as the continuous versions, so it was perhaps foolish of me to try to use the same names for both.
Nevertheless, I hope to show that there was method in my madness.

The objective now is to look at discrete and continuous versions of space(time) and momentum(-mass), in order to try and see how the various
concepts relate to each other, and if possible derive a discrete version of the Dirac equation. To do this 
properly requires moving up to the $24$-dimensional
representation on all of the fermions at once, in order to implement the  
continuous groups as left-multiplications and the finite group as right-multiplications, as explained in \cite{octahedral}.

The $24$-dimensional representation is the representation on a space spanned by the $48$ group elements, subject to the identification
of the group element $-1$ with the scalar $-1$. The easiest way to study it, therefore, is to study the action of the group on itself
by multiplication (on the left or right, as appropriate).
The continuous groups are then obtained by taking real linear combinations of the finite group elements.
Since the individual finite group elements represent individual particle interactions, the linear combinations represent
suitable weighted averages of particle interactions to describe the environment. In this way the model implements a form
of Rovelli's relational quantum mechanics (RQM).

Before attempting to work with the full $24$ dimensions, let us see how this works in some smaller representations.
First look at the complex $2$-space spanned by $P$ and $R$, with the Dirac matrices corresponding to $i,j,k,d$ as listed in
Equation~(\ref{DiracPR}). In this case it is easy to see that $i+d$ is a matrix with a single non-zero entry, and deduce that the
algebra generated by $i,j,k,d$ is the full matrix algebra $M_2(\mathbb C)$. A similar calculation shows that the algebra generated
by the matrices given in Equation~(\ref{DiracHI}), corresponding to $i,j,k,d$ acting on $H,I$, also generate the full matrix
algebra $M_2(\mathbb C)$. 

In both cases, therefore, we can construct a group $SL(2,\mathbb C)$, which contains $i,j,k$ but not $d$.
Indeed, $i,j,k$ lie inside $SU(2)$, and
we can construct a group $U(2)$ containing $i,j,k$ and $d$. In the standard model, the finite group is not used, so that the
groups $SL(2,\mathbb C)$ and $U(2)$ appear to be independent. 
Nevertheless, in the standard implementation, $SL(2,\mathbb C)$ and the scalar $U(1)$ act continuously,
while $SU(2)$ is reduced essentially to the single element $id$, representing $\gamma_5$.

To see how this works, note that in the finite group algebra model $SL(2,\mathbb C)$ acts on two copies of the same 
type of spinor, while in the standard model it acts on two complex conjugate spinors. This can be accommodated by changing
the sign of $\sqrt{-2}$ on one of the spinors, at the cost of introducing some extra complications. First of all, it converts the
scalar $\sqrt{-2}$ into an element that acts as $\sqrt{-2}$ on half the spinor and as $-\sqrt{-2}$ on the other. In standard
language, it swaps the scalar $i$ with $i\gamma_5$, a feature that 
occurs in other
extensions of the standard model \cite{Clifford}. 

Secondly, this implies that the standard model uses
the eigenspaces of $id$ to try to define these two halves of the spinor. But the eigenspaces of $id$ define quite a different copy
of $U(1)$, so that the standard model then has to explain the difference between the two. This may be where many of the
complications in the implementation of electro-weak mixing come from.

\subsection{Mass}
\label{mass}
As we have seen, the complex structure on $P,Q,R,S$ defines a generation of electrons. The given complex structure splits
the spinor into eigenspaces of $id$, one spanned by $R$ and $P-Q$, the other spanned by $S$ and $P+Q$.
There is therefore some discrete vector in one of these spaces that defines the electron. It is not obvious which is the best
vector to use, 
and it probably depends on which of the triplet symmetries $w, iw, jw,kw$ we use for the electron generation symmetry.

If $w$ represents colour symmetry, then we should probably use $w^i=jw$ for electron generation symmetry, although
$iw$ and $kw$ are also possible. Now as a permutation
\begin{eqnarray}
jw&=&(P,S,-R).
\end{eqnarray}
The fixed point, $Q$, should probably represent
the common feature of the three generations, namely the charge. 
Also, there should be components in
both eigenspaces, representing the left-handed and right-handed spinors. 

This analysis suggests $\pm R-(P+Q)$ and $\pm S+(P-Q)$ as likely candidates.
I would suggest that the mass should be encoded in $P,R,S$ somehow, with the largest contribution probably from $P$, and
therefore suggest taking the case in which the coefficients of $P$ and $S$ are both positive.
In this case, the three generations would be (up to permutations) 
 \begin{eqnarray}
 e&=&-Q-R+S,\cr 
 \mu&=&-Q+P-R,\cr\tau&=&-Q+S+P.
 \end{eqnarray}
Hence a proton might be
represented by $P+Q-R+S$, assuming it is independent of electron generation. Then we must surely take the neutron, as the
neutral version of the proton, as $P-R+S$. Adding together the spinors for
the three generations of electrons and three protons then gives the spinor 
\begin{eqnarray}
e+\mu+\tau+3p&=&5(P-R+S)\cr
&=& 5n.
\end{eqnarray}
Hence we see a possible justification for the mass 
equation (\ref{emutau}), including a
reason for appearance of the number $5$ there.

The mass itself must be defined by some 
real linear combination of $P,Q,R,S$, whose 
physical meanings we still have to elucidate.
But the given equation is independent of what the actual mass values are. It is therefore a clue to some
quantum structure of matter that is more fundamental than mass. At this stage, we are still working
with a toy model, but the equation is sufficiently well justified by the model that it can 
be claimed as a
prediction of the mass of the $\tau$ particle, with a relative standard uncertainty of $2\times 10^{-8}$, compared to a current experimental uncertainty
\cite{perspective,CODATA}
of $9\times 10^{-5}$:
\begin{eqnarray}
m(\tau)_P &=&  1776.84145(3) \mbox{ MeV}/c^2,\cr
m(\tau)_E &=& 1776.82(16) \mbox{ MeV}/c^2.
\end{eqnarray}

\subsection{The environment}
\label{environment}
In the full model, the generation symmetry $w$ also acts on $P,Q,R,S$, so that the matrix algebra generated by $i,w,d$ acts as the full
matrix algebra $M_4(\mathbb R)$, and the copy of $SL(2,\mathbb C)$ used in the standard model has no particular significance.
By the principles of RQM, it encodes certain macroscopic properties of the environment. Our task, therefore, is to try to identify
which properties of the environment these are, and how they are encoded. 

Since the ambient group is now $SL(4,\mathbb R)$, or $GL(4,\mathbb R)$,
we must expect that these environmental properties are gravitational in nature. There is no well-defined copy of either $SL(2,\mathbb C)$
or $SO(3,1)$, and
therefore no well-defined concept of inertial frame in this model. There is, however, a well-defined copy of $SO(4)$, in case we wish to use it.

In the full model, the $4\times 4$ matrices model the tensor product of a Euclidean spacetime with its dual.
The dual here is the Euclidean equivalent of Lorentzian $4$-momentum, which I take to be mass-momentum.
Now we have a choice of interpretations, according to whether we quantise the rows or the columns or both or neither.
Classical mechanics quantises neither, and quantum mechanics in principle quantises one and not the other, but in
practice is somewhat ambivalent. The finite group model is at its most useful also when one is quantised and not the other,
and enforces a rather more rigorous approach to separating the discrete from the continuous, but allows us some
choice in interpreting which is which. 

If we want a model that is close to the standard model, we probably want left-multiplication by $SL(4,\mathbb R)$ to act on
columns representing mass-momentum. Under certain circumstances, in which mass can be treated as constant,
one can then change the interpretation to
columns of $4$-momentum. Then we must have the finite group acting by right-multiplication on rows representing
quantised Euclidean spacetime. These rows also  
constitute one half of the Dirac spinor, representing an
elementary particle of some kind. In other words, quantised spacetime in this model is, essentially, part of the Dirac spinor, with its
complex structure removed.

\section{Modelling of physical processes}
\subsection{Lithium}
The mass equation (\ref{emutau}) involves three leptons, and the smallest `classical' object that contains three leptons is a lithium atom.
Now there is a problem with lithium in the standard model, because the amount of lithium that should theoretically be produced in 
Big Bang nucleosynthesis is a factor of $3$ or $4$ larger than that which is inferred from observations of very old stars. This is known as the cosmological lithium problem, and a number of different hypotheses have been put forward as possible explanations \cite{lithiumproblem}.

If the given equation is not just a curious coincidence, but has a real physical meaning, then it suggests that there may be some dynamical process
by which one side of the equation can be converted into the other, not just as an abstract mass, but as a physical collection of particles. Such a
process does not change the mass, and is not electromagnetic, so in quantum terms it must be a strong force reaction. But it is not a strong
force reaction in the standard model, because it fails to preserve lepton or baryon number. It would therefore have to be a ninth dimension
of the strong force extending the gauge group from $SU(3)$ to $U(3)$.

Extensions of the standard model 
by adding a new $U(1)$ gauge group
have been considered before in the literature, for example in the theory of 
sphalerons \cite{sphaleron}. These however
preserve the difference $B-L$ between baryon number $B$ and lepton number $L$, while the proposed interaction does not preserve $B-L$,
unless we add in five neutrinos.
So a sphaleron does not seem to meet the criteria. What we require is a ninth massless gluon. Such a gluon is the gauge boson of an extension to the
strong force, with a gauge group $U(1)$ that is definitely different both from electromagnetic $U(1)$ and from $U(1)_Y$.

It should be noted that although in my model both sides of the equation have $15$ spins, 
in the standard model the left hand side has six and the right hand side has five. To translate to the standard model it 
may be possible to add two neutrinos to the left hand side and three anti-neutrinos to the right hand side
in order to balance the spins. However, one would need to test this reaction experimentally in order to
determine whether this is a reasonable interpretation.
At the moment all we can do is note the difference and hope that it can be resolved
eventually, by experiment rather than theory. A model of the type I propose contains this single interaction
that is impossible in the standard model, but seems otherwise capable of
reproducing the entire standard model. We must therefore look at predictions of consequences of this interaction,
in order to compare with experiment.

What happens in the theory, however, is that the spins measured with respect to internal space, external space and parametric space
are not necessarily the same thing. The $3\times 3$ matrices that describe the $9$ gluons in internal space are the same matrices that
describe the $3$ dimensions of momentum for $3$ generations of neutrino in external space. The translation converts (spin) angular momentum to
momentum, and therefore acts as a super-symmetry to convert between bosons and fermions. But this super-symmetry, unlike others that
have been proposed, acts only on massless particles that travel at the speed of light, so between gluons and neutrinos.
Moreover, this supersymmetry is not a physical process, but a purely mathematical operation, that translates between two different
ways of looking at the same thing. In other words, in my model gluons and neutrinos are the same thing.
Moreover, the ninth gluon has both an interpretation as a graviton (in internal space) and a combination of three neutrinos
travelling in perpendicular directions (in external space).

Now we
can see why this gauge boson 
is essentially undetectable 
in existing experiments. It represents an interaction that can occur in the Big Bang, 
but that cannot be detected in the Large Hadron Collider. 
If we want it to explain lithium disappearance, then presumably 
it has to convert the leptons and protons into neutrons, not the other way round. Therefore it has to take place in an environment where there is enough energy and enough neutrinos to bump two electrons up to higher generations simultaneously,
and at the same time enough pressure to squeeze them close enough to the protons to get in range of the strong force. That is quite a tall order, and explains
why it has never been seen in experiments.

This is not the place to try and quantify this proposed reaction. There is no prediction that can be tested in the laboratory, and there is a great
deal of uncertainty in the astronomical consequences. But if there is a probability given by a 
relationship between two copies of $U(1)$ in $U(3)$, expressed by the
CP-violating phase, then it may not be completely irrelevant that the cosine of this angle is approximately $1/3$.
If this is the probability that a lithium atom arises at a particular stage in Big Bang nucleosynthesis, when the standard model
assumes this probability is $1$, then we may have gone some way towards solving the cosmological lithium problem.

One further remark is that if this process can be used in one direction to reduce the amount of lithium produced in the Big Bang,
then it can in principle run in reverse to increase the amount of lithium produced in stars. For this one 
would apparently need an element with
at least five neutrons in its nucleus, such as beryllium-9. If these five neutrons convert to $e+\mu+\tau+3p$, and these particles 
are then shaken apart violently enough by thermal motion, then the weak force could perhaps take hold and split the nucleus apart.
If the muon and tau particle decay to electrons, then an enormous amount of energy is released, as two entire nucleons have effectively
been turned into pure energy.

By combining this process with an inverse beta decay, it might also be possible to ignite it in a lithium-7 atom. Then the process would run as
\begin{eqnarray}
\label{neutrondestruction}
4n &\rightarrow& \mu+\tau+2p,
\end{eqnarray}
with some energy input/output from neutrinos and/or antineutrinos, 
after which the $\mu$ and $\tau$ would decay to electrons with the release of energy and more neutrinos and anti-neutrinos.
Initiating this reaction might only require an energy input 
of around $1$ MeV, but it might be very rare, and therefore not detected, because of the low interaction rates of neutrinos.
A process of this kind could in principle be happening inside the Sun's corona, which contains significant amounts of highly ionised metals
at very high temperatures, and a huge flux of neutrinos. 
If so, then the process releases so much energy that it could maintain the temperature of the corona at its
observed value in the millions of degrees. In other words, this is a process that might solve the coronal heating problem
\cite{corona}. 

Indeed, there is one plausible candidate explanation in terms of re-connecting
magnetic field lines, which 
might indeed be explained at a quantum level by the process I have described.
The explosive conversion of neutrons into electrons and protons opens up magnetic field lines that can then connect
to the pre-existing magnetic field lines.
The reason why the process can only occur (on a large scale) in the corona is that it is mediated by the graviton.
It therefore needs the lithium to be gaseous and completely ionized in order to remove all the electromagnetic forces
that would otherwise prevent it from happening. The same process might also 
take place
in heavier elements that are gaseous and sufficiently highly ionised.

If this is true, then in any large collection of lithium atoms one might just by chance see this reaction taking place occasionally even at
low average temperatures. Even a relatively small number of such reactions could release enough energy to cause a lithium battery to
overheat and catch fire.
While the previous suggestions in this section cannot be tested experimentally, this one can. All one needs to do is bombard lithium atoms
with neutrinos or anti-neutrinos, tuned to the right energy, and see what happens. But 
one should be prepared for an explosion.

\subsection{Heavy metal}
In order to maintain an equilibrium between $e+\mu+\tau+3p$ and $5n$ in the nucleus, an atom needs a protective shell of electrons
to keep the nuclear leptons from escaping. It is clear that if the number of copies of $5n$ is greater than the number of copies of $3p$, then no such equilibrium can exist. Therefore the ratio of neutron number to proton number must be less than $5/3\approx 1.67$. This is known to be
the case, and can be explained by the standard model even without the graviton. But it is possible that the proposed graviton may
simplify the explanation, since it permits the development of a model that
does not require the input of a large number of experimentally determined parameters.

The neutron/proton number ratio is typically around $1$ for small atoms, and increases to around $1.59$ for uranium.
But there is significant variation from this overall trend. This means that different chemical elements respond differently to the graviton.
Now almost all of this difference is attributed in the standard model to the strong force, so there is very little left that could in principle
be detected by purely gravitational experiments. Nevertheless, the model proposed in this paper implies that there is a residual
difference between the gravitational attractions between different chemical elements that might be detectable by experiments that
measure Newton's gravitational constant $G$ using different materials.

To estimate the sensitivity of experiment that would be required, let us assume that an anomaly arises from the mass difference between
$n$ and $e+p$, and compare iron and lead as being close to the two extremes of typical materials that might be used. The ratio of atomic weight to
atomic number is around $2.148$ for iron, and $2.527$ for lead. Normalizing iron to the lead scale would give it an effective atomic number of $22.1$,
thereby increasing its effective mass by approximately $3.9$ copies of $n-e-p$, that is around $.0032$ amu, equivalent to $58$ ppm.

Now CODATA has never claimed an accuracy as good as this before 2014, when a relative standard uncertainty
of $46$ ppm was given. More recent experiments claim to reach $12$ ppm, but are in conflict with each other
\cite{bigG1,bigG2}.
What my model proposes is a possible source of systematic error that has not been acounted for.
Moreover,  predicted errors from this source are of the correct order of magnitude to account for the observed discrepancies.
Or rather, the model actually proposes a mixing of the strong force with gravity, which implies that $G$ is not in reality a
constant. An experiment to compare the relative gravitational forces due to diamond (weight/number ratio 2.00) and platinum
(ratio 2.50), to test a predicted difference of $76$ ppm in the 
values of $G$,
would certainly be 
eye-catching! 

This raises the question of how accurate a theory of gravity (such as Newtonian gravity or general relativity)
 that is based on an absolute constant $G$ can actually be. Astronomical evidence suggests that it works well on a quite large distance scale, 
 and on a scale 
 of gravitational accelerations down to very low accelerations, but that it starts to break down
 on the scale of the acceleration of the Solar System towards the centre of the Milky Way.
 This must be the scale, therefore, on which the quantum effects of the graviton must be taken into account,
 and the differences between internal, external and parametric spacetimes must be taken into account also.
 In other words, we must take account of the finite speed of propagation. 
 
 This analysis suggests that the gravitational anomalies that result in the hypothesis of dark matter may not in fact be physical
 anomalies at all, but artefacts of our theory. To resolve these anomalies, we must understand what the variations in $G$ that we
 appear to detect on Earth are actually telling us about a mixing, or coupling, between the theories of the strong force and gravity
 that have up to this point developed largely independently. Since the model explains this coupling in terms of a supersymmetry
 between gluons and neutrinos, we must look for clues in the PMNS matrix. My guess would be that this matrix is related to the
 geometry of the relationship between our orbit around the Sun, and the Sun's orbit around the centre of the galaxy, but I have
 no direct evidence to justify this conjecture at this point.
 What may be happening is a mixing between the neutrinos emanating from the Sun, and the neutrinos emanating from
 the rest of the galaxy, and it is this mixing that the experiments detect.

\subsection{Noether's theorem}
One of the fundamental results in classical physics is Noether's Theorem, that derives a conservation law from a symmetry group.
But it applies only to continuous field theories, so that we need a different kind of result to obtain conservation laws in a discrete theory.
It so happens that I proved a suitable result in representation theory in my Ph.D. thesis, and published it a few years later \cite{1.2.3}.

Translated into physical language, it states that in any model of the type I have been considering,
any quantity that is globally invariant is locally conserved. 
Quite apart from any mathematical proof, this principle would seem to be of fundamental physical
significance in any attempt to unify theories of the very large and the very small. 
It should be noted that the concepts of `local invariance' and `global conservation' are logically inconsistent, even if sometimes
of practical utility, so cannot be used in a unified theory. 
Also, note that the principle does \emph{not} say that a locally conserved quantity is necessarily globally invariant.

This principle applies, for example, to electric charge, which is globally invariant in all theories of electrodynamics, and is experimentally
confirmed to be conserved in all known interactions. 
On the other hand, 
the `generation' of an elementary fermion is known not to be conserved in weak interactions, which implies that it is not globally
invariant. 

This property has also been confirmed in experiments on neutrinos. That is, the generation of a given neutrino is indeed
measured differently by different observers in many circumstances. What is different about my model is that this phenomenon
is an immediate consequence of Noether's theorem, and quite independent of any properties of mass. 
Moreover, it applies also to the generations of electrons, and therefore has the astonishing consequence that mass is not
invariant in this model.

The immediate reaction may therefore be to reject this model as impossible, but given its successes so far that might be
a rather rash response. Rather, we should examine the experimental evidence for and against this conclusion. It is clear that
mass is invariant in the theory of special relativity, and therefore conserved in quantum electrodynamics. The assumption that mass
is invariant in general relativity is a form of the equivalence principle, which has often been questioned theoretically, 
although direct experimental evidence against it is lacking. 

Nevertheless, indirect experimental evidence for the failure of the equivalence principle
is abundant in astronomical data, for example in the rotation curves of galaxies. In the previous section
I have suggested ways in which direct experimental evidence may now be within reach.
In a later section I also present some circumstantial evidence that the measured masses of elementary
particles may depend in some subtle ways, not yet understood, on properties of the gravitational field.

Presumably the physical process by which a neutrino changes its generation can be modelled as an interaction with
the proposed graviton, which changes spin properties but not mass properties. There are two other experimentally
observed processes of this kind that I have drawn attention to elsewhere \cite{octahedral}. One is the oscillations of neutral
kaons between three observable states, without any directly detectable change in mass
\cite{CP,kaonanomaly,Kaon2}. Another is a potential influence of
the gravitational field on the muon gyromagnetic ratio
\cite{muong-2,muontheory,muonHVP,Fodor}. In all cases there seems to be a primary dependence on the direction
of the gravitational field, but there may also be a secondary dependence on latitude and longitude separately.

\section{Further remarks on representations}
\label{remarks}
\subsection{Complex spacetime}
\label{Cspacetime}
Now the real representation of $GL(2,3)$ on $A,B,C,D$ is equivalent to that on $P,Q,R,S$,
although the integral representations are quite different. Hence, if a complex structure on $P,Q,R,S$ is
necessary for defining the Dirac spinor for an electron, then this same complex structure must also be defined on $A,B,C,D$.
In other words, the choice of the first-generation electron leads to an associated complexification of quantum spacetime,
although it is not at all clear what physical interpretation we should give to such a complexification. 

In this case, the action of $d$ on $A,B,C,D$ is the negative of its action on $P,Q,R,S$, which changes the sign of $\sqrt{-2}$,
giving
\begin{eqnarray}
\sqrt{-2}A = C-D,&&
\sqrt{-2}C=B-A,\cr
\sqrt{-2}B= -C-D,&&
\sqrt{-2}D=A+B.
\end{eqnarray}
However, it is not clear how many of these choices are just conventions. 

In this context, $A,B,C,D$ are generators for quantum spacetime, corresponding to $W,X,Y,Z$, in which $W$ has been identified
as being related to mass, energy and/or time in some way. As generators for the integral representation, it is not enough to
take $A,B$ and their multiples by $\sqrt{-2}$, but we need either $A,C$ or $A,D$. If we take $A,C$, then the imaginary part
consists of $C-D$ and $B-A$. On a quantum scale, these differences are discrete versions of imaginary derivatives, 
$i\hbar d/dt$ and $i\hbar d/dx$, for some suitable direction $x$. In particular, $C-D$ corresponds to the energy term
$i\hbar d/dt$ for the first generation of electron. If so, then $B-C$ and $D-B$ correspond to the energy terms for the
other two generations.

Now the question is, what is the relationship, if any, between the quantum `space' directions $B,C,D$, and macroscopic
space directions? The three generations of electron seem to be related 
to the directions $B-C$, $C-D$ and $D-B$,
that 
span a $2$-space, not a $3$-space. Then it might be possible to take $B+C+D$ as some arbitrary
direction in space, so that the three generations do not depend on direction, but on a combination of directions.

\subsection{Modulo $3$}
\label{GL23mod3}
The real purpose of this exercise is to prepare for the reduction modulo $3$, where it is no longer possible to
keep the spinors separate from the spacetime representations. It is here that we hope to see some link between the spinors and
spacetime that might throw some light on the $3$-generation structure of fermionic matter.
Now we can calculate the following relations modulo $3$:
\begin{eqnarray}
H+J&\equiv& A-C+D \pmod 3,\cr
-H+I+J-K&\equiv& B+C+D\pmod 3,\cr
I-K&\equiv& -P+R-S \pmod 3,\cr
H+J-I-K&\equiv& Q+R+S \pmod 3.
\end{eqnarray}

From this we obtain a splitting of the spinor into `left' and `right' components, whose complicated
relationship to 
the standard model concepts of left and right has already been discussed. 
This splitting does not exist in the real or integral
representations, but can be expressed in terms of our definition of $\sqrt{-2}$ as follows:
\begin{eqnarray}
(1+\sqrt{-2})H &\equiv & -(1+\sqrt{-2})P\pmod 3,\cr
(1+\sqrt{-2})I &\equiv& (1+\sqrt{-2})(P+Q)\pmod 3,\cr
(1-\sqrt{-2})I &\equiv& (1-\sqrt{-2})A\pmod 3,\cr
(1-\sqrt{-2})(I-H)&\equiv& (1-\sqrt{-2})B \pmod 3.
\end{eqnarray}
In other words, the two components of the $H,I,J,K$ spinor can be thought of as the multiples by $1\pm\sqrt{-2}$.
Let us therefore call the multiples of $1+\sqrt{-2}$ \emph{positive}, and the multiples of $1-\sqrt{-2}$ \emph{negative},
to avoid conflict with the left/right distinction(s) in the standard model.
We then have a correlation between the positive parts of the $H$-spinors and the $P$-spinors, while the negative
part of the $H$-spinor correlates with quantum spacetime. The negative part of the $P$-spinor, however, has no such
correlation with any other part of the representation.

Now dividing through by $1\pm\sqrt{-2}$, and using the fact that $(1+\sqrt{-2})(1-\sqrt{-2})=3$, we have
\begin{eqnarray}
H&\equiv & -P \pmod {1-\sqrt{-2}},\cr
I&\equiv& P+Q\pmod {1-\sqrt{-2}},\cr
H&\equiv& A-B \pmod {1+\sqrt{-2}},\cr
I&\equiv& A \pmod {1+\sqrt{-2}}
\end{eqnarray}
The first two of these congruences express a gluing between the $H,I,J,K$ spinor and the $P,Q,R,S$ spinor, while the
last two seem to express a gluing between the $H,I,J,K$ spinor and spacetime.

\subsection{The PIMs modulo $3$}
\label{pims}
The reduction of the given $12$-dimensional representation modulo $3$ can be written as the direct sum of two PIMs in nine
different ways. Altogether this representation has $28$ distinct submodules, of which $14$ are indecomposable, so that 
some computational assistance is recommended at this point. The `Meat-axe' package \cite{Meataxe,Ringe},  which can also
be found incorporated into GAP \cite{GAP} and MAGMA \cite{MAGMA}, was
designed specifically for computing with modular representations of finite groups, and is ideally suited for this purpose.

It turns out that there is just one PIM that contains $A,B,C,D$, and the remaining two generators can be taken as
\begin{eqnarray}
E=\begin{array}{|c|ccc|} \hline.&.&-&-\cr.&+&.&-\cr.&+&-&.\cr\hline\end{array},&&
F=\begin{array}{|c|ccc|} \hline.&-&-&.\cr.&.&-&+\cr.&-&.&+\cr\hline\end{array}.
\end{eqnarray}
Similarly, there is just one PIM that contains $P,Q,R,S$, with two further generators

\begin{eqnarray}
T=\begin{array}{|c|ccc|} \hline+&.&-&.\cr+&-&.&.\cr.&.&+&+\cr\hline\end{array},&&
U=\begin{array}{|c|ccc|} \hline.&-&.&.\cr+&.&.&-\cr+&+&+&.\cr\hline\end{array}.
\end{eqnarray}
In each case, these vectors must be interpreted only over $\mathbb F_3$, and no meaning is attached to any 
particular pre-image over $\mathbb Z$.

\subsection{The Green ring}
\label{Greenring}
To understand the representation theory of $GL(2,3)$ over $\mathbb F_3$
 in enough detail for the applications to physics, we need to understand the tensor products
of all indecomposable representations. This information
is encapsulated in the Green ring \cite{Greenring}.
The representations $3^\pm$ and the full-height towers are projective, and the tensor product of
a projective with anything is always projective, so that these tensor products can be determined from the Brauer character table.
So all we really need to know about tensor products is what is left over after the projectives have been thrown away. 

In the case we have here, we need only two pieces of information besides the Brauer character table, namely
\begin{eqnarray}
2^+\otimes 2^- \rightarrow 1^+, && \begin{matrix}1^+\cr 1^-\end{matrix} \otimes \begin{matrix}1^-\cr 1^+\end{matrix} \rightarrow 1^+.
\end{eqnarray}
All other required information can be obtained by tensoring with irreducibles, and using the characters.
In the above cases, for example, the character table tells us that the actual tensor products are
\begin{eqnarray}
2^+\otimes 2^- = 1^++3^-, &&\begin{matrix}1^+\cr 1^-\end{matrix} \otimes \begin{matrix}1^-\cr 1^+\end{matrix} = 1^+
+ \begin{matrix} 1^-\cr 1^+ \cr 1^-\end{matrix}.
\end{eqnarray}

The cited paper by Benson and Parker \cite{Greenring} contains a wealth of information about the Green ring, analysed in great detail,
and includes in particular a deep analysis of the structure of the `glue' that holds the indecomposable representations together.
The `glue' in this sense is the difference between an arbitrary representation and the direct sum of its irreducible constituents.
They give as examples, among others, the groups $Sym(2)$, $Sym(3)$ and $V_4$ that are important for our purposes, for all
primes dividing the group order.

It is likely that these Green rings will repay further study, 
in particular for the understanding, or at least modelling, of how gluons behave in particle interactions.
Unfortunately, although the Green ring for $GL(2,3)$ over $\mathbb F_3$ is quite 
straightforward, that over $\mathbb F_2$ is much more
complicated, as there are infinitely many indecomposables.
This fact is likely to prove an obstacle to
implementing a comprehensive theory of electro-weak interactions on this model.

\section{
Experimental evidence}
\subsection{The nature of mass}
So far, the model has enough detail to make one specific theoretical prediction, that is the formula (\ref{emutau}),
in addition to a proposed extension to the strong force, from $8$ gluons to $9$,
that leads to speculative explanations of 
a number  of unsolved problems.
With current experimental values for the five masses involved, the two sides of the equation agree to
an accuracy of 6ppm, compared to a relative standard uncertainty in the data of 26ppm. 
But since the right hand side is known to much higher accuracy than the left hand side, this permits a
prediction of the mass of the tau particle to three more significant figures than are currently known. 

With more work to assign
specific particles to specific elements of the model, it should be possible to make more predictions of this kind.
Although I have done the calculations for this particular equation inside an integral version of the real
representation $4^0$ (in the notation of \cite{octahedral}) for simplicity, they
should really be done in the representation $4^0\otimes 4^0$, which splits into symmetric and anti-symmetric parts
\begin{eqnarray}
S^2(4^0) &\cong & 1^++3^++3^++3^-,\cr
\Lambda^2(4^0) &\cong& 1^-+2^0+3^-.
\end{eqnarray}
Thus I am really working in the anti-symmetric square, ignoring the $2^0$ component. As an integral representation,
rather than a real representation, there is no direct sum decomposition of this form, but there is a $4$-dimensional
quotient module consisting of the mass representation $3^-$ glued on top of the charge representation $1^-$, obtained
by ignoring the colours in $2^0$. 

To find other equations of this type, we should expect to require the full $6$ dimensions.
For example, there should be an analogous equation that relates the three boson masses (Higgs, $W$ and $Z$) to the scalar
neutron mass. We would like to find such an equation in 
an integral version of $1^-+2^0+3^-$, that is 
equivalent to a monomial representation
on $6$ points. Some trial and error leads to the following suggestion for vectors to represent the particles:
\begin{eqnarray}
Z^0&=&(0,1,1;0,1,1)\cr
W^+&=&(1,0,1;1,0,1)\cr
W^-&=&(-1,1,0;1,1,0)\cr
H^0&=&(0,1,1;0,0,0)\cr
n&=&(0,0,0;1,1,1)
\end{eqnarray}
so that
\begin{eqnarray}
\label{Higgsmass}
Z^0+W^++W^-&=&(0,2,2;2,2,2)\cr
&=& 2H^0+2n.
\end{eqnarray}
Of course, this is pure guesswork at this stage, and may be completely wrong,
but this equation is at least consistent with experiment. Moreover, it provides a prediction for the mass of the Higgs boson
that is more tightly constrained than the current experimental value. 
More speculatively, the process (\ref{neutrondestruction}), when combined with weak force processes, can be
represented as the annihilation of two neutrons, so that (\ref{Higgsmass})
suggests it may be possible to regard this process as being mediated by
the Higgs boson. However, it is not clear that this is a helpful way to interpret these equations.

More generally, the model suggests that there should be approximately four fundamental
particle masses, from which the $15$ fundamental masses in the standard model can be derived by simple linear equations.
Looking for such equations in advance of developing a theory is dangerous, and can easily degenerate into numerology.
On the other hand, if some sufficiently convincing equations could be found, they might
usefully guide the development of the theory. The equation(s) found so far do not in fact relate the standard model
fundamental masses to each other, but also involve the proton and the neutron. This begs the question of whether we might
sensibly allow in some other particles as well. But the more particles we consider, the greater the danger of falling into
the trap of numerology.

\subsection{Universal mass equations}
The best defence against numerology is to restrict to a small number of particles, whose mass is known as precisely as possible.
We may therefore wish to exclude the quarks (especially the light quarks), and replace them by mesons or baryons. One possibility
would be to add in the remaining six members of the baryon octet, in which we already know one more equation, namely the Coleman--Glashow
relation. Another possibility would be to take the six masses of the up/down/strange
pseudoscalar mesons.
Either case would give us a representation $2^0+3^++3^-$, so would extend into the symmetric square. To get both copies of $3^+$
we may need both the mesons and the baryons. Or perhaps the Higgs boson equation should really be expressed in terms of $3^++3^-$
inside the symmetric square, or more plausibly in a mixed copy of $1^++2^0+3^+$,
in which case we might get away with using only the mesons. 

In this case, we can restrict our attention to the 
$12$ distinct masses of the following particles:
\begin{eqnarray}
e,\mu,\tau,n,Z,W^\pm,\pi^0,\pi^\pm,K^0,K^\pm,\eta,\eta'.
\end{eqnarray}
For present purposes, the $3$ neutrino masses are too small to be relevant, so that 
these $12$ masses are a suitable replacement for the $15$ fundamental masses in the standard model.
Now we have three tasks: 
\begin{enumerate}
\item to express the six quark masses in terms of the $12$ chosen masses;
\item
to find approximately $8$ simple equations relating these $12$ masses;
\item to find approximately $4$ equations relating some of these masses to properties of the gravitational field and/or acceleration
of the experiment and/or observer.
\end{enumerate}

The first of these tasks was considered in \cite{FQXI}, where the following suggestions were presented, in terms of
the quark masses $u,d,s,c,b,t$:
\begin{eqnarray*}
e+u&=&d\cr
\mu&=&s+2d\cr
\tau&=& c+5s\cr
c&=&s+K^++\eta+\pi^0\cr
b+s+d+\pi^++2\pi^0&=&5n\cr
t+c+d+\eta'&=& b+s+u+Z^0+W^+.
\end{eqnarray*}
In order to provide some motivation for these equations, note that the first three attempt to relate the generation
symmetry of the electrons to the up/down/strange symmetry of quarks, but it turns out that the symmetry is
mostly shifted to down/strange/charm symmetry, with a residual up/down/strange symmetry mixed in.
In terms of representation theory, this seems to be relating two copies of $3^-$, one for electrons and one for quarks,
together with $1^+$ and $1^-$ for mass and charge.

The fourth equation relates the charm and strange quarks via three pseudoscalar mesons,
one of each type in the meson octet. Here again the meson triplet relates to an up/down/strange triplet,
but in a different guise, so perhaps represented in $3^+$. It should be noted that there is an alternative version of this equation in which the positive charge lies in the pion rather than the kaon.
The fifth equation expresses the 
sum of the three generations of down quark, in terms of three pions and the scalar neutron mass. Again we are relating
two natural triplets, presumably represented in $3^+$ and/or $3^-$.
Finally, the sixth involves all six quarks: the more massive of
each weak doublet is on the left hand side, and the less massive on the right hand side, and the relationship
involves the weak force $Z$ and $W$ bosons. Here we see two triplets, a singlet and a doublet, so the representation 
is probably $1^++2^0+3^++3^-$.

One can easily check that these equations can be solved uniquely for $u,d,s,c,b,t$. 
From readily available data one can then 
calculate the values of the
masses of the two sides of these equations, which are here tabulated in units of MeV/$c^2$:
\begin{eqnarray}
\begin{array}{ll}
2.8\pm .7 \pm .5 & 4.8\pm .5\pm .3\cr
106 & 105\pm 5\cr
1777 & 1750\pm 35\cr
1275\pm 25 & 1271\pm 5\cr
4690\pm 30 & 4698\cr
175450\pm510\pm710 & 175844\pm 30
\end{array}
\end{eqnarray}
Clearly the first equation is very dubious, and of no real value. The potential value of the others as clues to further development of the model
is also unknown: the proof of the pudding is in the eating. They are however consistent with experiment, and have some predictive
value in terms of obtaining more accurate values for the masses of the three heavy quarks (charm, bottom and top). 

\subsection{Evidence for non-invariance of mass}
The second task was not addressed in \cite{FQXI}, but six suggestions were made for the third task, which might in principle
enable us to derive one or two suggestions for the second task. Here we need to consider what properties of the gravitational field
and/or acceleration we are prepared to use. As a general principle, we should only use dimensionless parameters, and we should only
use parameters that are constant to high accuracy over large distances, and that can be measured accurately. However, these principles
are very difficult to satisfy simultaneously in practice. The choice that was made in \cite{FQXI} was not to compromise on dimensions or
accuracy, but to compromise a little on constancy. This surely implies that the suggested equations are wrong. I present them
in the hope that they may provide clues to assist in finding the correct equations.

The four dimensionless almost-constant parameters that were chosen were two temporal ratios, namely the ratios of the periods of revolution of the 
Earth around the Sun, and the Moon around the Earth, to the period of rotation of the Earth on its own axis; and two spatial angles, namely the angles of
tilt of the last two motions with respect to the first. These parameters all affect the acceleration of the observer relative to the experiment, and
although the effects would seem to be negligible on 
the laboratory scale, there is no guarantee that they 
are 
negligible on a quantum scale, since
there is no guarantee that they 
are 
continuous on  
that scale.

For the sake of argument, let us take the two frequency ratios to be 365.24 and 29.53, and the two angles to be $23.44^\circ$ and $5.14^\circ$.
The six mass ratios (out of only $11$ independent mass ratios available)
that were considered have the following experimental values:
\begin{eqnarray}
n/p&\approx& 1.00137842\cr
e/p&\approx&5.44617021\times 10^{-4}\cr
\pi^\pm/\pi^0 &\approx& 1.03403\cr
K^\pm/K^0 &\approx& .99202\cr
W^\pm/Z^0 &\approx& .88146\cr
K^0/\eta&\approx& .90893
\end{eqnarray}

If we are looking for circumstantial evidence that these mass ratios may not be globally invariant, then what we are looking for is ways
in which these numbers look similar to simple functions of the chosen four dimensionless parameters of the gravitational field.
In this context, it is hard to avoid noticing the similarity of the mass ratios to the numbers
\begin{eqnarray}
1+1/(2\times 365.24) &\approx& 1.0013689\cr
\sin(23.44^\circ)/(2\times 365.24) &\approx&
5.44543 \times 10^{-4}\cr
1+1/29.53&\approx& 1.03386\cr
\cos^2(5.14^\circ)&\approx& .99197\cr
\cos(23.44+5.14)^\circ &\approx& .87815\cr
\cos(23.44^\circ)\cos(5.14^\circ) &\approx& .9137
\end{eqnarray}
although the last two are certainly less convincing than the first four.

Clearly, some or all of these numerical coincidences may indeed be coincidences with no relevance to physics.
But the foundations for a model that I have presented in this paper certainly imply that some equations of this
general type must hold. 
These algebraic foundations are not sufficient to 
explain any of these equations in particular, 
and 
none of these equations can be exact, although 
the level of approximation
should be broadly consistent with what the above figures suggest.
To what extent these suggestions support the main thesis of this paper is a moot point, but if no such coincidences could be found,
the model would have been in trouble. 

\subsection{Black holes and cosmic expansion}
\label{blackholes}
It goes without saying that a `black hole' in the strict sense of a singularity cannot exist in any quantum theory. But in the sense of a very large amount of matter concentrated in a very small amount of space, a black hole is not impossible in quantum theory. On the other hand, since no signals can get out of a black hole, it is impossible for us to observe what is going on inside. Nevertheless, we must assume that the laws of physics 
inside a black hole are exactly the same as the laws we observe in the rest of the universe.

In this paper I have proposed only one addition to the known laws of physics, and that is a coupling of gravity to quantum physics via a graviton
that mediates a specified quantum process. This process is capable of explaining in qualitative terms a number of mysteries ranging from neutrino
oscillations and the properties of kaons to coronal heating and the `missing' primordial lithium. Let us therefore investigate what happens
if we add these gravitons to a black hole.

For this purpose it does not much matter how we identify the graviton with anything in the standard model of particle physics, only that
it is a massless `particle' that travels at the speed of light. The proposed process by which a graviton can be detected in a very strong
gravitational field, for example in the interior or on the surface of a neutron star, converts the energy of the graviton into the energy of
two neutrinos and two anti-neutrinos. Instead of the spin $2$ graviton of a quantum gravity related to general relativity, 
what we
see is a set of four entangled spin $1/2$ particles. For practical purposes, 
we might as well work with the neutrinos, which
we know to exist, rather than with the 
graviton, which may or may not exist as a particle in its own right. The `graviton' in this
sense is a highly non-local concept, as its essential properties reside in the entanglement of the neutrinos.

Then, since the gravitons 
convey information, they cannot travel faster than light, so they cannot escape from the black hole. That means that the
normal process of exchange of gravitons is disrupted, and the gravitons only flow one way. Presumably this weakens gravity outside the black hole,
and strengthens it inside the black hole. 
But there is a subtlety here, due to the non-locality of the information. If the neutrinos fall into the black hole, it may still be possible, or even
mandatory, for the
anti-neutrinos to escape, or vice versa. 

Moreover, this is a dynamic process that drives the evolution of an entire galaxy over its entire lifetime.
This implies that the gravity we measure in the Solar System at a distance of 27000 light-years from the centre of a galaxy that is
more than 13 billion years old depends on those parameters in an essential way.
It is therefore somewhat remarkable that a theory of gravity that does not mention those parameters describes the motion of much larger
structures than the Solar System so well. What is not remarkable is that it fails to account for the motion of structures in which the
acceleration becomes less than the acceleration of the Solar System towards the 
galactic centre \cite{Milgrom1,Milgrom2,Milgrom3,MOND,MOND2,newparadigm}.

If gravity increases within the black hole, then the collapse proceeds at an increasing rate, until the central core becomes essentially a
quark-gluon plasma, but now with $9$ gluons instead of $8$. The ninth gluon, moreover, converts quarks to leptons, which causes 
the core to explode. For this process to occur, we must assume that the total number of quarks and gluons in the plasma is very large, and
therefore the explosion is very large. It is like the Big Bang, but without the inflationary stage. It must be assumed, therefore, that it
produces a result very much like the result of the Big Bang, but on a smaller scale.

The model predicts that every large galaxy, and perhaps even the small ones, ends up like this, so that 
it might in principle be possible to build a model of the universe in which the Big Bang is replaced by a myriad of
Small Bangs. However, it would be quite challenging for such a model to explain \emph{all} the cosmological observations
that the Big Bang Theory can explain. 
On the other hand, a Small Bang Theory could explain cosmic voids, and it would not need an inflationary stage, so could be built
entirely on known physics, together with one additional process that corresponds to a ninth gluon.

Moreover, if cosmic expansion is driven by these Small Bangs pushing against each other, then it could explain an expansion that is not uniform,
and not constant. This might throw some light on the Hubble Tension
\cite{Hubble1,Hubble2,nodarkenergy}, for example, and apparently anomalous redshifts of some objects
in the universe. In particular, it may be worthwhile to look in detail at whether a Small Bang Theory can explain the high relative velocities
of the galaxies within the Virgo cluster, that is revealed by a wide spread of redshifts above and below the large-scale average.
There are of course many other aspects of the formation of the universe we see that would be different if the driving force were
Small Bangs taking place within an existing universe, rather than one Big Bang outside of time and space. It may be, of course, 
that the Small Bang Theory 
can be easily ruled out by some such aspects, such as properties of the cosmic microwave background.
However, initial investigations look promising, and suggest that more work could profitably be done on this model.

\subsection{Entanglement}
In the previous section I proposed that the conventional assumption of a spin $2$ graviton be replaced by an entangled set of
two neutrinos and two anti-neutrinos, in such a way that the gravitational properties are contained in the entanglement, rather than in the
individual particles. If this is a valid interpretation, then it suggests a more general coupling between entanglement and gravity.
In other words, we could turn the picture around, and instead of using entanglement to explain quantum gravity, we could
use quantum gravity as an explanation of entanglement.

The bizarre properties of entangled electrons or photons \cite{Rae}, 
fully confirmed by experiment, appear to rule out any local theory of quantum mechanics,
so that non-locality must be built into the theory somewhere. My proposal is to use quantum gravity as a source for this non-locality, physically
realised in terms of interactions with neutrinos. Since even photons are known to interact with the gravitational field, this is certainly
physically possible, although it is hard to see how it could ever be experimentally tested. Nevertheless, my suggestion is that 
entanglement is maintained only as long as the gravitational field remains constant.

\section{Conclusion}
\label{conclusion}
In this paper I have looked at various possibilities for quantising elementary particles in terms of 
finite and integral representations of finite groups. The integral representation theory is in principle 
what would be required for a complete theory, intermediate between the continuous (real or complex) representation theory
that is
required for quantum mechanics and the standard model as they currently exist, and the finite (modular) representation
theory, localised at the primes $2$ and $3$, that seems to be
required for the underlying discrete structure of spin doublets, weak doublets, colour triplets and generation triplets. 

However, integral representation theory is very hard, even for very small groups, so that I have concentrated on the
modular representations in order to restrict attention to the conceptual problems, rather than the numerical problems.
I have been mainly concerned with the $3$-modular representation theory, as the $2$-modular theory is a significant
source of the difficulties in the integral theory.
The primes $2$ and $3$ appear both as (multiplicative) orders of elements in the group, and as (additive) orders of
elements in the underlying field, so that there are doublets and triplets of both kinds. 

There is therefore no
possibility of separating the doublets into one (weak) theory and the triplets into another (strong) theory. The two are inextricably
linked, in a variety of different ways, both in the group itself, and in the coefficients that are used in the representations. 
The prime $2$ in the group
is associated with weak doublets, and in the field is associated with left-handed and right-handed spinors,
so that although these pairs of concepts are closely related, 
they are mathematically not the same thing,
as indeed is also clear from physical experiment.
The prime $3$ in the group is associated with generation triplets, and in the field is associated with colour triplets,
so that again the mathematics 
explains why these two sorts of triplets are
physically very different.  

The $3$-modular representations, acted on by the elements of order $2$, explain why the weak force is chiral,
but do not explain which chirality it has. The $2$-modular representations, acted on by elements of order $3$,
explain why there are three generations of fermions, but do not explain the 
mass hierarchy. The $3$-modular
representations, acted on by elements of order $3$, explain colour confinement.
 The $2$-modular representations, acted on by elements of order $2$, explain 
why left-handed and right-handed spinors are different. But to combine all these ingredients into a unified theory requires
a deep understanding of the integral representation theory of $GL(2,3)$, and has not been attempted here.

My intention is rather to draw attention to a part of mathematics that may be of use in such an endeavour, but which
is not well-known to physicists. Over a period of more than ten years, and a series of several papers \cite{perspective,Clifford,
finite,octahedral}
I have scoured the field of algebra for useful ingredients for potential unified theories, 
and have rejected almost everything. This includes all the parts of algebra
that have been prominent in many previous attempts to construct unified theories of various kinds, including unitary groups such as the
Pati--Salam \cite{PatiSalam} four-colour
$SU(4)$, the Georgi--Glashow \cite{GG} $SU(5)$ and larger examples of unitary groups
\cite{HH},  orthogonal groups and associated Clifford algebras
\cite{Furey1,Furey2,Furey3},
exceptional Lie groups and associated Lie algebras and Jordan algebras \cite{E6,E8,DG}, quaternion and octonion algebras
\cite{Dixon}, finite simple groups
such as the Monster and its associated vertex operator algebra \cite{Borcherds,FLM}, 
this and other types of `moonshine', other sporadic groups such as the
Conway group and the Leech lattice, and a host of other speculative ideas. 

All of these ideas have shown promise
at some time or other, and some are still under active investigation. But my conclusion after all this work 
is that the part of algebra that shows the most
promise for genuinely useful  applications to fundamental physics is the representation theory, real, complex, integral and modular, of the group
$GL(2,3)$. There is, of course, no guarantee that a viable theory can be built on this foundation. But it appears to be the only part of
algebra that both has a reasonable chance of success and has not already been exhaustively explored in the physics literature.
It is therefore, I suggest, worthy of serious consideration by those physicists who are interested in going beyond the standard model.

To support this conclusion I have explored a number of potential applications of the model to well-known problems in physics,
including the cosmological lithium problem, the coronal heating problem, neutrino oscillations and kaon decays.
While these applications are currently speculative, they all derive from
a single proposed new reaction, mediated by a ninth, colourless, scalar gluon. The fact that the model can explain 
several
apparently disparate problems by the same physical process is a strong argument in its favour.

\end{document}